\documentclass[10pt]{amsart}
\usepackage{amsmath,amsthm, amssymb,amsfonts, amscd, epsfig, epsf,euscript} %showkeys

\theoremstyle{plain}
\newtheorem{thm}{Theorem}[section]
\newtheorem{lem}[thm]{Lemma}
\newtheorem{pro}[thm]{Proposition}
\newtheorem{cor}[thm]{Corollary}
\theoremstyle{definition}
\newtheorem{Def}[thm]{Definition}
\newtheorem{rem}[thm]{Remark}
\newtheorem{exa}[thm]{Example}
\numberwithin{equation}{section}
\newcommand{\mc}[1]{\mathcal{#1}}
\newcommand{\mbf}[1]{\mathbf{#1}}
\newcommand{\mfr}[1]{\mathfrak{#1}}
\newcommand{\mr}[1]{\mathrm{#1}}

\newcommand{\coh}[3]{ H ^{#1} \left (  #2, \,  #3\right )}

\newcommand{\lclcoh}[3]{H ^{#1}_{#2} (   #3 )}
\newcommand{\erpq}[3]{E_{#1}^{#2, #3}}
\newcommand{\ext}[3]{\mathrm{Ext}^{#1}(#2, \, #3 )}
\newcommand{\homdot}[2]{\mathrm{Hom}^{\bullet}(#1, \, #2 )}

\newcommand{\ceil}[1]{\left\lceil{1}\right\rceil}

\DeclareMathOperator{\lcm}{lcm} 
\DeclareMathOperator{\reg}{reg}

\DeclareMathOperator{\Spec}{Spec}

\newcommand{\sign}{{\mathrm{sign}}}
\DeclareMathOperator{\vregnum}{reg-num_{\bf v}}
\DeclareMathOperator{\vjregnum}{reg-num_{{\bf v}_j}}
\DeclareMathOperator{\vIregnum}{reg-num_{{\bf v}_I}}
\newcommand{\dreg}{{\mathrm{dreg}}}

\newcommand{\Ass}{{\mathrm{Ass} }}

\newcommand{\V}{{\mathbb V}}
\newcommand{\Z}{{\mathbb Z}}
\newcommand{\R}{{\mathbb R}}
\newcommand{\C}{{\mathbb C}}
\newcommand{\F}{{\mathbb F}}
\newcommand{\N}{{\mathbb N}}

\newcommand{\proj}[1]{{\mathbf P}^{#1}}

\begin{document}
\title[Regularity and resolutions]{Regularity and resolutions for multigraded modules}
\author{Jerome W. Hoffman and HaoHao Wang}
\address{Department of Mathematics\\
              Louisiana State University\\
              Baton Rouge, Louisiana 70803 \\}
\email{hoffman@math.lsu.edu}
\address{    Department of Mathematics \\
              Southeast Missouri State University \\
              Cape Girardeau, MO 63701\\}
\email{hwang@semo.edu}

\begin{abstract}
This paper is concerned with the relationships between two
concepts, vanishing of cohomology groups and the structure of free
resolutions. In particular, we study the connection between
vanishing theorems for the local cohomology of multigraded modules
and the structure of their free multigraded resolutions.
\end{abstract}


\thanks{We thank Jessica Sidman for her
comments and suggestions. The second author thanks NSF-AWM and
GRFC for travel support and the hospitality of Louisiana State
University while visiting the first author.}

\maketitle

\section{Introduction}
\label{S:intro}
Castelnuovo-Mumford regularity was first defined by Mumford
\cite{M} for coherent sheaves $\mc{F}$ on projective spaces
$\proj{n}$. $\mc{F}$ is $m$-regular if \[H^i(\proj{n},
\mc{F}(j))=0, \   \ \forall \  \  i+j \geq m, \   \ \forall  \ \ i
\geq 1.\] Mumford showed that the regularity of a sheaf measures
the smallest degree twist for which the sheaf is generated by its
global sections. More generally, regularity of a sheaf has
implications for the structure of a resolution of that sheaf by
direct sums of the sheaves $\mc{O}(j)$.

There is a corresponding notion for graded modules.  Regularity of
$\Z$-graded modules was investigated by many people, notably,
Bayer-Mumford \cite{BM}, Bayer-Stillman \cite{BS}, Eisenbud-Goto
\cite{EG}, and Ooishi \cite{O}. Let $S=K[x_1, \dots, x_n]$ be the
polynomial algebra in $n$ variables over a field $K$ with the
standard grading.  If $M$ is a finitely generated graded
$S$-module, $B = (x_1, ..., x_n)$ is the maximal ideal, then $M$
is $m$-regular if
\[H^i_B(M)_j=0,  \   \   \forall \  \  i+j \geq m+1.\]
An important result is:
\begin{thm}(\cite{BM})\label{T:BM}
Suppose $K$ is a field and $M$ is a graded $S$-module. Then $M$ is
$m$-regular if and only if the minimal free graded resolution of
$M$ has the form
\[0 \to  \displaystyle{\bigoplus _{\alpha = 1}^{r_{s}}}S(-d_{\alpha
,s}) \to  \cdots  \to \displaystyle{\bigoplus _{\alpha =
1}^{r_{1}}} S (-d_{\alpha ,1}) \to \displaystyle{\bigoplus_{\alpha
= 1}^{r_{0}}}S(-d_{\alpha ,0}) \to M \to 0 \] where $d_{\alpha ,i}
\le m + i$ for all $i \ge 0$.
\end{thm}
This result gives an explicit relationship between the regularity
of a graded module, in other words the vanishing of certain graded
pieces of local cohomology $H^i_{B}(M)$, and degree bounds on
generators of the free graded resolution of $M$.

Recently, there has been a great interest in studying regularity
of multigraded modules over a multigraded algebra. Geometrically
such algebras arise as the coordinate algebras of projective toric
varieties. Let $S$ be a $G$-graded $K$-algebra where $K$ is a
field and $G$ is a finitely generated abelian group. Let $B$ be a
monomial ideal in $S$ and $M$ a finitely generated $G$-graded
$S$-module. Multigraded regularity relates the vanishing of
$H^i_B(M)_{\mbf{d}}$ for $\mbf{d}$ in certain convex regions of
$G$ with the structure of resolutions of $M$. In case $(S, B) $ is
the homogeneous coordinate ring and irrelevant ideal of a smooth
projective toric variety $X$, $G$ is the divisor class group of
$X$, and the geometric version of the theory imposes vanishing
conditions on $\coh{i}{X}{\mc{F}\otimes \mc{O}(\mbf{d})}$. A
special case of this, with $G = \Z ^2$, was done by the authors
\cite{HW}, and the general case is due to Maclagan and Smith
\cite{MS}. Here we introduce the definition of multigraded
regularity as in \cite{MS}.

Let $\mc{C} =\{\mbf{c}_1, \dots, \mbf{c}_l\}$ with $\mbf{c}_i \in
G$ be a fixed subset of $G$, and let $\N \mc{C}$ be the
sub-semigroup of $G$ generated by $\mc{C}$
($\N = \Z _{\ge 0}$).
Define
\begin{equation}
\label{E:nc} \N \mc{C}[j] = \bigcup_{\mbf{w}=(w_1, \dots, w_l) \in
\N^{l}, \  \ \sum w_i = |j| } (\sign(j)( w_1{\mbf{c}_1}
+\dots+w_l{\mbf{c}_l} ) +\N \mathcal{C}).
\end{equation}
Then $M$ is $\mbf{m}$-regular with respect to $B$ and $\mc{C}$, if
for each $i \geq 0$:
\begin{equation}
\label{E:reg} H^i_B(M)_{\mbf{d}}=0,  \  \ \forall \mbf{d} \in
\mbf{m}+ \N\mc{C}[1-i],
\end{equation}
and we denote this by $\mbf{m} \in \reg_{B, \mc{C}}(M)$. There is
no loss of generality in assuming that $S = K [x_1, ..., x_n]$ is
a free algebra. The degrees of the generators $x_i$ are elements
$\mbf{a}_i \in G$ and we let   $\overline{\mbf{a}}_i$ be their
images in $G \otimes _{\Z}\R$. The set $\overline{\mc{A}}=
\{\overline{\mbf{a}}_1,..., \overline{\mbf{a}}_n  \} $ is an
integral vector configuration in $G \otimes _{\Z}\R$. One assumes
\begin{itemize}
\item[1.] The monomial ideal $B$ corresponds to a maximal cell
$\Gamma \subset G \otimes _{\Z}\R$ in the vector configuration
$\overline{\mc{A}}$ via
\[B :=  \langle \prod _{i\in \sigma } x_i : \sigma \subseteq [n] \mr{\ with\ }
  \Gamma \subseteq \mr{pos}( \overline{\mbf{a}}_i : i \in \sigma \rangle,
  \ [n] := \{ 1, ..., n\} \rangle
\]
where the positive hull of a set of vectors is
\[ \mr{pos}( \mbf{v}_1, ..., \mbf{v}_m )
:= \{ \lambda _1 \mbf{v}_1+ ... +\lambda_m \mbf{v}_m: \lambda _i
\in \R_{\ge 0}.  \}\] The choice of this maximal cell $\Gamma$ is
implicit in all the results of \cite{MS}. \item[2.]
$\mr{pos}(\overline{\mc{A}})$ is a pointed cone with
$\overline{\mbf{a}}_i \ne \mbf{0}$ for all $i$.
\end{itemize}
These two conditions are satisfied by $(S, B)$, the homogeneous
coordinate ring and irrelevant ideal of a projective toric
variety. For the main results connecting to the structure of
resolutions one assumes
\begin{itemize}
\item[3.] $\mc{C} \subseteq \mc{K}$ where $\mc{K} $ is the Kaehler
cone (for detailed definition see \cite[Section 2.]{MS}).
\end{itemize}
Among their results are certain bounds on the generators and
syzygies of an $\mbf{m}$-regular module, but these are not as
sharp as those of the classical theory. In particular they obtain,
under the further hypothesis,
\begin{itemize}
\item[4.] $\N \mc{C}= \mc{K} = \mc{K}^{\mr{sat}}$
\end{itemize}
not a free resolution, but a complex with $B$-torsion homology
augmenting to the module $M$ (see \cite[Theorem 7.8]{MS}). Since
modules with $B$-torsion homology become zero when converted to
sheaves on the corresponding toric variety, this is adequate for
applications for sheaf theory. But it is still of interest to
obtain more precise results about the structure of free
resolutions. In \cite{HW} a concept of ``strong regularity'' was
introduced which gave results analogous to Theorem \ref{T:BM}. One
of the themes of this paper is to study the effect of stronger
regularity conditions on the structure of free resolutions.

We first investigate the relationship between certain regions
\[\mc{D}^i_B(J) \subset \mc{Z}^i_B(M) =
 \{\mbf{d} \in G : H^i_B(M)_{\mbf{d}}=0\}
\]
and the sets $J =\{J_0, J_1, \dots, J_s\}$ in a minimal free
resolution
\[ 0 \to \bigoplus_{\mbf{d}_s \in J_s} S(-\mbf{d}_s) \to \cdots \to
\bigoplus_{\mbf{d}_0 \in J_0} S(-\mbf{d}_0)\to M \to 0 .
\]
of a $G$-graded $S$-module $M$. This question has two aspects. On
the one hand, starting with $J$, one would like to find easy to
describe regions $\mc{D}^i _B(J)$ with the property that any
module $M$ with a resolution of this type with degrees in $J$ has
$\lclcoh{i}{B}{M}_{\mbf{d}}= 0$ for $\mbf{d} \in\mc{D}^i _B(J) $.
On the other hand, starting from known regions $\mc{D}^i \subset
G$ of vanishing of $\lclcoh{i}{B}{M}_{\mbf{d}}$ one would like to
find a free resolution for $M$ with the $p$-th syzygy module
having degrees in regions $J_p$ which can be described explicitly
in terms of the regions $\mc{D}^i$.  We will discuss this
relationship in Section \ref{S:relation} after we introduce some
necessary tools in Section \ref{S:tool}.

In practice one would like the regions $\mc{D}^i_B(J)$ to be the
intersection of $G$ with a rational convex polyhedral region
in $G\otimes \R$. The faces of these would be defined by linear
inequalities. Thus
one is led to consider regularity with respect to regions defined
by $h \ge 0$ where $h$ is a linear form. In effect this means that
one is coarsening the multigrading to a simple, but generally
nonstandard grading. This has been considered in a recent paper by
Sidman-Van Tuyl-Wang \cite{SVW} who extended some of the
techniques used in \cite{HW}.  They introduced a notion of
positive coarsening vector which coarsens the multidegrees on $M$
to a single nonstandard $\Z$ grading, and obtain finite bounds on
the multidegrees of a minimal free resolution of $M$. If $\phi:  G
\to G_0$ is a surjection, we define a $G_0$-grading of $M$ by
\[ M_{\mbf{d}_0} = \bigoplus_{\mbf{d} \in \phi^{-1} (\mbf{d}_0)}
M_{\mbf{d}}.\] One considers vanishing regions contained in $G_0$
for $H^i_{B_0}(M)$ with respect to a $G_0$-graded ideal $B_0$, and
a set $\mc{C}_0 \subseteq G_0$.  The simplest case is when
$G_0=\Z$, and if $G=\Z^r$, $\phi(\mbf{d})=\mbf{v}\cdot \mbf{d}$
for a coarsening vector $\mbf{v} \in G$. This case is studied in
detailed in \cite{SVW}.  Often, the $(B, \mc{C})$-regular region
can be estimated in terms of $(B_0, \mc{C}_0)$-regularity regions,
defined by various $\phi$'s, and the advantage of considering the
$\phi$'s is that one can often give easily described constraints
on the regions $J_p$. We will discuss this in Section
\ref{S:vdot}.

Another approach is to consider regularity with respect to a
family of $(B_0, \mc{C}_0)$ defined by natural geometric
conditions.  One example we consider is the canonical
decomposition of the nullcone:
\[Z= \V(B)= \bigcup_{i=1}^t Z_i,  \  \ Z_i = \V(B_i),\]
where $(S, B)$ is the homogeneous coordinate ring and irrelevant
ideal of a simplicial toric variety $\mbf{P}_{\Sigma}$ defined by
the fan $\Sigma$.  The ideals $B_i$ are linear.  We employ a
Mayer-Vietoris spectral sequence in local cohomology to relate
vanishing regions of $H^i_B(M)$ to those of $H^i_{B_I}(M)$ with
respect to ideals $B_I$, where $B_I$ is defined by the
intersection of the various $Z_i$. This method works especially
well for fans of the form $\Sigma= \Sigma_1 \times \cdots \times
\Sigma_s$. In Section \ref{S:family}, we will discuss this issue,
and also we illustrate its use for multiprojective spaces
$\proj{m_1} \times \cdots \times \proj{m_t}$.

The example of multiprojective spaces shows the necessity of
considering gradings $\phi$ which assign $\mbf{0}$ to some of the
variables, which was not treated in \cite{SVW}.  We can include
such gradings by writing $S=R[x_j]$ where $R=K[x_i]$ where $x_i$
are the variables with $\phi(\deg(x_i))=0$, and $x_j$ are the
variables that are $\phi$-positive, i.e., admissible for the
theory of regularity.  For this, we need to generalize the theory
developed in \cite{MS} to $G$-graded algebras over a commutative
Noetherian ring $R$, not just a field.  We do this in Section
\ref{S:0degree} following ideas of Ooishi \cite{O}. We end this
paper by giving some illustrative examples in Section
\ref{S:examples}.

\section{Tools}
\label{S:tool} This section recalls some well-known spectral
sequences and derives some corollaries needed later.

Let $\mc{A}$ be an abelian category. We consider chain complexes,
indexed by $\Z$, over $\mc{A}$ both of homological type:
\[\begin{CD}
C_{\bullet}= \cdots @>>>C_n @>\partial>> C_{n -1}@>>> \cdots
\end{CD},
\]
and of cohomological type:
\[\begin{CD}
C^{\bullet}= \cdots @>>>C^n @>d>> C^{n +1}@>>> \cdots
\end{CD}.
\]
These are formally equivalent by defining $C^n = C_{-n}$. We let
$C[p]^{\bullet} $ be the shifted complex:
 $C[p]^{n} = C^{p+n}$, and we often identify
an object $M$ in $\mc{A}$ with a complex concentrated in degree
$0$.

Given complexes $A^{\bullet}$, $B^{\bullet}$, the first bounded
above and the second bounded below, we recall that the groups
$\ext{i}{A^{\bullet}}{B^{\bullet}}$, which exist if $\mc{A}$ has
enough injective objects. This is computed by taking a
quasi-isomorphism $B^{\bullet}\cong I^{\bullet} $ where
$I^{\bullet}$ is a bounded below complex of injective objects.
Then
\[\ext{i}{A^{\bullet}}{B^{\bullet}} =
    H^i (\homdot{A^{\bullet}}{I^{\bullet}})
\]
where $\homdot {A^{\bullet}}{I^{\bullet}}$ is the total hom
functor defined in \cite[pp. 63-64]{HRD}.

\begin{pro}
\label{P:ss1} There is a spectral sequence
\[\erpq{1}{p}{q} = \ext{q}{A^{\bullet}}{B^{p}} \Longrightarrow
\ext{p+q}{A^{\bullet}}{B^{\bullet}}.
\]
This converges if $A^{\bullet}$ is bounded above and $B^{\bullet}$
is bounded below.
\end{pro}

\begin{proof}
This is the hypercohomology spectral sequence of the functor $T=
\homdot{A^{\bullet}}{-}$ applied to the complex $K = B^{\bullet}$
with the filtration $F = \sigma _{\ge p}$, whose graded terms are
$B^{p}[p]$, \cite[pp. 20-21]{DHII}.
\end{proof}

\begin{cor}
\label{C:ss1} If
\[\begin{CD}
\cdots @>>> C_2 @>>> C_1 @>>>C_0 @>\varepsilon>> M @>>> 0
\end{CD}
\]
is a resolution, then there is a convergent spectral sequence
\[
\erpq{1}{-p}{q} = \ext{q}{A^{\bullet}}{C_{p}} \Longrightarrow
\ext{q-p}{A^{\bullet}}{M}.
\]
\end{cor}

\begin{proof}
By assumption, the module $M$, regarded as a complex concentrated
in degree $0$, is quasi-isomorphic with the complex $C_{\bullet}$,
and we may thus replace $M$ by that complex in the computation of
Ext. Apply the proposition and shift to upper numbering in
indexing: $p \to -p$.
\end{proof}

\begin{pro}
\label{P:ss2} There is a spectral sequence
\[\erpq{1}{-p}{q} = \ext{q}{A^{p}}{B^{\bullet}} \Longrightarrow
\ext{q-p}{A^{\bullet}}{B^{\bullet}}.\] This converges if
$A^{\bullet}$ is bounded above and $B^{\bullet}$ is bounded below.
\end{pro}

\begin{proof}
We let $K^{\bullet} = \homdot{A^{\bullet}}{I^{\bullet}}$. Since
$I^{\bullet}$ is a complex of injective objects, the functor
$C^{\bullet} \to\homdot{C^{\bullet}}{I^{\bullet}} $ is exact. Thus
we can define a decreasing filtration by
\[F^p K^{\bullet} = \homdot{A^{\bullet}/\sigma _{\ge
1-p}A^{\bullet}}{I^{\bullet}}.
\]
Moreover,
\[\mr{Gr}_{F}^{p} K^{\bullet} =
\homdot{\mr{Gr}_{\sigma}^{-p}A^{\bullet}}{I^{\bullet}} =
\homdot{A^{-p}[p]}{I^{\bullet}} = \homdot{A^{-p}}{I^{\bullet}}[-p]
\]
The spectral sequence of a filtered complex
\[\erpq{1}{-p}{q}= \mr{H}^{q-p}(\mr{Gr}_{F}^{-p}K^{\bullet})
\Longrightarrow \mr{H}^{q-p}(K^{\bullet})
\]
gives our spectral sequence, in view of
\begin{align*}
\mr{H}^{q-p}(\mr{Gr}_{F}^{-p}K^{\bullet}) &=
\mr{H}^{q-p}(\homdot{A^{p}}{I^{\bullet}}[p]\\
&= \mr{H}^{q} (\homdot{A^{p}}{I^{\bullet}})=
\ext{q}{A^{p}}{B^{\bullet}}.
\end{align*}
\end{proof}

Let $X_*$ be a simplicial topological space, that is, a
contravariant functor from the category $\Delta$ whose objects are
the totally ordered sets $[n] = \{ 1, ..., n\}$ and whose
morphisms preserve the order.  An augmented simplicial topological
space $a: X_{*}\to S$ is of cohomological descent if the canonical
map of functors from $\mr{D}^+ (S) \to \mr{D}^+ (S) $ :
\[\varphi : Id \to \mr{R}a_{\ast}a^{\ast}\]
is an isomorphism (see \cite[\S 5]{DHIII}). We will be interested
in the following special case: Let $j:Z\hookrightarrow X$ be the
inclusion of a closed subspace, and
\[Z = \bigcup _{i=1}^{t}Z_i\]
is a union of closed subsets. This defines an augmented simplicial
space $a: W_{*}\to Z$ by taking the usual Cech complex:
\[W_k = \coprod _{\# I = k+1} W_I, \quad W_I = \bigcap _{i \in
I}Z_i, \quad I \subseteq \{1, \ldots, t  \}\]  Since $\coprod Z_i
\to Z$ is a proper surjective map, $a$ is of universal
cohomological descent (\cite[\S 5]{DHIII}). Let $a_k : W_k \to Z$
be the canonical map. The components of this, $a_I : W_I \to Z$
for $\#I = k+1$, are closed imbeddings. Thus, $a_{k\ast}$ is an
exact functor on the category of abelian sheaves, and hence
$\mr{R}^i a_{k\ast} = 0$, for $i> 0$. The means that, for any
abelian sheaf $F$ on $Z$,
\[F \longrightarrow C^{\ast}(W_{*}, F)\]
is a resolution, where
\[C^{k}(W_{*}, F) = a_{k\ast}a_{k}^{\ast} F = \bigoplus _{\# I =
k+1}a_{I\ast}a_{I}^{\ast} F\] with the usual Cech differential
(\cite[\S 5]{DHIII}). Since $j_{\ast}$ is an exact functor on
abelian sheaves, we obtain a resolution
\[j_{\ast}F \longrightarrow j_{\ast}C^{\ast}(W_{*}, F)
\quad\text{with}\quad j_{\ast}C^{k}(W_{*}, F) = \bigoplus _{\# I =
k+1}F_{W_{I}, X}\] where $F_{W_{I}, X}$ is the sheaf on $X$ that
coincides with $F \mid W_{I}$ on $W_{I}\subset X$ and is extended
by $0$ on the rest of $X$.

Apply this to the constant sheaf of integers $F = \Z _{Z}$. Recall
that, for any closed subset $Z\subset X$ we have a canonical
isomorphism:
\[
\lclcoh{i}{Z}{X, \, M} = \ext{i}{ \Z _{Z, X}}{M}
\]
where the left-hand side is local cohomology, and the right-hand
side is the Ext group for the category of abelian sheaves on $X$
(\cite[Expos\'e I, Prop. 2.3]{aG}).

\begin{pro}
\label{P:mayer-vietoris} Let $Z = \bigcup _{i=1}^{t}Z_i$ be a
union of closed subsets $Z_i \subset X$. Let $M$ be a sheaf of
abelian groups on $X$. There is a convergent spectral sequence
\[
\erpq{1}{-p}{q}= \bigoplus _{\#I=p+1} \lclcoh{q}{W_{I}}{X, \,
M}\Longrightarrow \lclcoh{q-p}{Z}{X, \, M}.\]
\end{pro}

\begin{proof}
The local cohomology is computed via the Ext's as explained. Since
$ \Z _{Z, X}$ is quasi-isomorphic to the complex
$j_{\ast}C^{*}(W_{*},\Z _{Z} )$ we may replace that sheaf by that
complex and apply  Proposition \ref{P:ss2}. We then obtain the
spectral sequence in the indicated form.
\end{proof}

When $t =2$ this spectral sequence reduces to the Mayer-Vietoris
sequence in local cohomology (\cite[p. 212, Ch. III, Ex. 2.4]{H}).

In our application, $X=\Spec(S)$ where $S$ is a finitely generated
$G$-graded $R$-algebra, with $R$ being a commutative Noetherian
ring, and $G$ is a finitely generated abelian group. Equivalently,
we can think of $X$ as an affine scheme with an action of the
group scheme $T=\Spec R[G]$. The most important case is when
$G=\Z^r$ is free, then $T\cong \mbf{G}^r_m$, where $\mbf{G}_m$ is
the multiplicative group. We are interested in $G$-graded
$S$-modules $M$. They define quasicoherent sheaves $\tilde{M}$ on
$X$ also with a $T$-action.

We will be considering
\begin{align*}
& Z=\V(B)=\bigcup_{i=1}^t \V(B_i)= \bigcup_{i=1}^t Z_i\\
&W_I = \bigcap_{i \in I} Z_i =\V(B_I), \  \ \emptyset \neq I
\subseteq \{1, \dots, t\}
\end{align*}
where $B, B_i, B_I \subset S$ are $G$-graded ideals in $S$, well
defined up to radical. We have
\[\lclcoh{i}{B}{M}=\lclcoh{i}{Z}{X, \tilde{M}},
\text{  and  }\lclcoh{i}{B_I}{M}= \lclcoh{i}{W_I}{X, \tilde{M}}, \
\ \emptyset \neq I \subseteq \{1, \dots, t\}.\] These local
cohomology modules are $G$-graded.

\begin{cor}\label{C:1}
Fix $i\geq 0$ and  $\mbf{d}\in G$. Suppose that, for all
$\emptyset \neq I \subseteq \{1,2, \dots, t\}$, we have
\[
\lclcoh{i+\#I-1}{B_I}{M}_{\mbf{d}}=0,  \text{ then }
\lclcoh{i}{B}{M}_{\mbf{d}}=0.
\]
\end{cor}

\begin{proof}
This is an application of the above discussion and proposition
\ref{P:mayer-vietoris}. The terms of the spectral sequence
contributing to $\lclcoh{i}{B}{M}_{\mbf{d}}$ are the
$\lclcoh{i+\#I-1}{B_I}{M}_{\mbf{d}}$.
\end{proof}

Returning to the geometric situation:
\begin{pro}
\label{P:BIM} Let $Z\subset X$ be a closed subset. Suppose that
$M$ is a sheaf of abelian groups on $X$ with a resolution of the
form
\[
 \begin{CD} 0@>\partial
>>C_d @>\partial>>...@>\partial>> C_1
@>\partial>>C_0 @>\varepsilon>> M @>>> 0.
\end{CD}
\]
There is a convergent spectral sequence
\[ \erpq{1}{-p}{q}=\lclcoh{q}{Z}{X, C_p} \Rightarrow
\lclcoh{q-p}{Z}{X,M} \]
\end{pro}

\begin{proof}
Local cohomology is computed via Ext modules. Then apply Corollary
\ref{C:ss1}.
\end{proof}

In our application, such sheaves and their resolutions will come
from $G$-graded $S$-modules:

\begin{cor}
\label{C:BIM} Let $S$ be a $G$-graded $R$-algebra. Suppose
$X=\Spec(S)$, $M$ a $G$-graded $S$-module and suppose we have a
resolution by $G$-graded $S$-modules:
\[ \begin{CD} 0@>\partial
>>C_d @>\partial>>...
@>\partial>>C_0 @>\varepsilon>> M @>>> 0.
\end{CD}
\]
Let $B \subset S$ be a $G$-graded ideal. Fix $i\geq 0$ and $\mbf{d}
\in G$. Suppose that for all $p\geq 0$, we have $
\lclcoh{i+p}{B}{C_p}_{\mbf{d}}=0,  \text{ then }
\lclcoh{i}{B}{M}_{\mbf{d}}=0 . $
\end{cor}

\begin{proof}
Translating the sheaf language to module language, we see that the
terms of the spectral sequence in the proposition contributing to
$\lclcoh{i}{B}{M}_{\mbf{d}}$ are the
$\lclcoh{i+p}{B}{C_p}_{\mbf{d}}$.
\end{proof}

\begin{rem}
If $\proj{}_{\Sigma}$ is a smooth projective toric variety
associated with the fan $\Sigma$, then Cox \cite{C} has associated
to $\proj{}_{\Sigma}$ a multigraded algebra $S$ and a monomial
ideal $B \subset S$. In that case, we take $X = \Spec (S)$, with
its $T$-action as above, and $Z = \mathbb{V}(B)$. This case will
be discussed in Section \ref{S:family}.
\end{rem}

\section{regularity regions and the degrees of syzygies}
\label{S:relation}
We assume that $S = K[x_1, ..., x_n]$ is a $G$-graded
$K$-algebra, where $K$ is a field,
with $\deg x_i =\mbf{a}_i \in G$.
Let $Q$ be the semigroup of $G$ generated
by the $\mbf{a}_i$ for $1\leq i \leq n$.  We say $S$ is
positively multigraded by $G$, if $G$ is torsion-free,
$\deg x_i \neq \mbf{0}$ for all
$i$ and $Q$ has no non-zero invertible elements. This
implies that each $S_{\mbf{a}}$ is finite dimensional
for each $\mbf{a} \in Q$ (\cite[theorem 8]{SM}).
Any finitely generated $G$-graded
$S$-module $M$ has a finite free graded resolution
each of whose terms is a finite direct sum of modules of the form
\[
S(\mbf{d}_1)\oplus \ldots  \oplus S(\mbf{d}_t)
\]
(\cite[prop. 8.18]{SM}). If $S$ is positively multigraded by $G$,
then there is a well-defined notion of minimal free resolution,
and any finitely generated $G$-graded module will have a unique up
to isomorphism minimal resolution of this form. This is shown in
the paragraph before proposition 8.18 in \cite{SM}.

Recall the definition of regularity given in the introduction,
Equation \ref{E:reg}.
\begin{Def}
\label{D:Lreg} Let $B\subset S$ be a $G$-graded ideal in a
$G$-graded ring $S=K[x_1, \cdots, x_n]$ and let $M$ be a
$G$-graded $S$-module. Define:
\begin{itemize}
\item[1.] $\mc{Z}_B^i(M)=\{ \mbf{p} \in G :
H^i_B(M)_{\mbf{p}}=0\}. $ \item[2.] $ \reg_{B, \mc{C}}^i(M) =\{
\mbf{m} \in G: H^i_B(M)_{\mbf{p}}=0, \   \ \forall \mbf{p} \in
\mbf{m} + \N \mc{C}[1-i] \}. $ \item[3.] $ \reg_{B,
\mc{C}}(M)=\bigcap_{i\geq 0} \reg_B^i(M) .$ \item[4.] $M$ is {\it
$\mbf{m}$-regular with respect to $B$ and $\mc{C}$} if $\mbf{m}\in
\reg_{B, \mc{C}}(M)$.
\end{itemize}
\end{Def}

\begin{rem}
\label{R:L} Note that $\reg_{B, \mc{C}}^i(M)$ is the largest set
$A\subset G$ such that
\[
A+\N\mc{C}[1-i] \subseteq \mc{Z}_B^i(M).
\]
Thus, $\reg_{B, \mc{C}}(M)$ is the largest set $A\subset G$ such
that \[ A+\N \mc{C}[1-i]\subseteq \mc{Z}_B^i (M),  \   \ \forall i
\geq 0.
\]
\end{rem}

Our first aim is to estimate  $\reg_{B, \mc{C}}^i(M)$ in terms of
$\reg_{B, \mc{C}}^i(S)$ and the structure of a free graded
resolution of $M$. The point of this is that the vanishing regions
for $H^i_B(S)$ with respect to monomial ideals $B$ are effectively
computable. See the paper of Eisenbud, Mustata and Stillman
\cite{EMS} and also \cite[Theorem 2.1]{Mu1} and \cite[Proposition
3.2]{MS}.
\begin{Def}
\label{D:typeJ} Suppose that $J = \{ J_0, ..., J_s\}$ is a list of
finite subsets $J_i \subset G$.  A resolution of a $G$-graded
$S$-module $M$ of the form
\[ 0 \to \bigoplus_{\mbf{d}_s \in J_s} S(-\mbf{d}_s) \to \cdots \to
\bigoplus_{\mbf{d}_0 \in J_0} S(-\mbf{d}_0)\to M \to 0
\]
is called a resolution of type $J$.
\end{Def}
\begin{Def}
\label{D:modvan} Suppose that $J = \{ J_0, ..., J_s\}$ is a list
of finite subsets $J_i \subset G$, and suppose that we are given regions
 $\mc{V}_B^{i}(S)\subset \mc{Z}_B^{i}(S)$, for
all $i \ge 0$. Define
\[
\mc{D}^i _B (J) = \bigcap_{p =  0}^s \bigcap_{\mbf{d}_{p} \in J_p}
(\mbf{d}_{p} + \mc{V}_B^{i+p}(S)).
\]

\end{Def}

\begin{rem} \label{R:regions}
In the above definition, the notation does not reflect the choice
of the sets $\mc{V}_B^{i}(S)$, which is implicit. We could take
$\mc{V}_B^{i}(S)=  \mc{Z}_B^{i}(S)$, for every $i$, but the reason
for allowing greater flexibility in choosing the sets
$\mc{V}_B^{i}(S)$ is that the complete vanishing region
$\mc{Z}_B^{i}(S)$ might have a complicated form; or in any case a
``nice'' region, for instance, rational convex polyhedral subset
$\mc{V}_B^{i}(S)\subset\mc{Z}_B^{i}(S)$ might be known or
relevant. By the theorems of Eisenbud, Mustata and Stillman
\cite{EMS} such regions can often be found. In the lemmas that
follow we let $\mc{V}_B^{i}(S)\subset\mc{Z}_B^{i}(S)$, $i\ge0$ be
any selection of regions and $ \mc{D}^i _B (J)$ the corresponding
sets.
\end{rem}
Since
\[
\lclcoh{q}{B}{C_p}_{\mbf{d}}=\bigoplus_{\mbf{d}_p \in J_p}
\lclcoh{q}{B}{S}_{\mbf{d}-\mbf{d}_{p}},
\]
it follows that
\begin{equation}
\label{E:vanres} \mc{Z}_B^i(C_p) =\bigcap_{\mbf{d}_{p} \in
J_p}(\mbf{d}_{p}+ \mc{Z}_B^i(S)),   \  \reg_B^i(C_p)
=\bigcap_{\mbf{d}_{p} \in J_p}(\mbf{d}_{p}+ \reg_B^i(S) )  .
\end{equation}
\begin{lem}\label{L:LBI}
Suppose that $M$ has a resolution of type $J$. Then for each fixed
$i \ge 0$ we have
\[
\mc{D}^i _B (J)\subseteq \mc{Z}_B^i(M).
\]
\end{lem}
\begin{proof}
For any fixed $p\ge 0$ it follows from Equation \ref{E:vanres}
that\[ \mbf{d} \in \bigcap_{\mbf{d}_p \in J_p} (\mbf{d}_{p} +
\mc{Z}_B^{i+p}(S)) \Rightarrow  H^{i+p}_B(C_p)_{\mbf{d}}=0.
\]
If $\mbf{d}\in \mc{D}^i _B (J) $, this will hold for every $p\ge
0$, by the definition of the set $ \mc{D}^i _B (J) $.  If we
assume that $M$ has a resolution of type $J$, then Corollary
\ref{C:BIM} implies
\[
H^i_B(M)_{\mbf{d}}=0
\]
and therefore $\mbf{d} \in \mc{Z}_B^i(M)$.
\end{proof}

\begin{lem}\label{L:mreg}
Let $J = \{ J_0, ..., J_s\}$ be a list of finite subsets $J_i
\subset G$. Let $\mbf{m}\in G$ and suppose that for a fixed  $i
\geq 0$,
\[ \mbf{m} + \N \mc{C}[1-i] \subseteq \mc{D}^i _B (J).\]
Let $M$ be a $G$-graded $S$-module with a resolution of type $J$.
Then $\mbf{m}\in \reg_{B, \mc{C}}^i(M)$. If this holds for every
$i \ge 0$, then $\mbf{m}\in \reg_{B, \mc{C}}(M)$.
\end{lem}

\begin{proof}
The result follows by Lemma \ref{L:LBI} and the definition of
regularity \ref{D:Lreg}.
\end{proof}

The next series of lemmas estimate the regularity of $M$ in terms
of the regularity of $S$.
\begin{lem}\label{L:regJ}
Let $J = \{ J_0, ..., J_s\}$ be a list of finite subsets $J_i
\subset G$. Let $\mbf{m} \in G$. Suppose that $M$ is a $G$-graded
$S$-module with a resolution of type $J$. If that for a fixed $i
\geq 0$,
\[\mbf{m} + \N \mc{C}[1-i] \subseteq \bigcap_{p =0}^s
\bigcap_{\mbf{d}_{p} \in J_p} (\mbf{d}_{p} + \reg^{i+p}_{B,
\mc{C}}(S)+\N \mc{C}[1-i-p]),\] then $\mbf{m}\in \reg_{B,
\mc{C}}^i(M)$. If this holds for every $i \ge 0$, then $\mbf{m}\in
\reg_{B, \mc{C}}(M)$.
\end{lem}

\begin{proof}
 By Remark \ref{R:L}
\[
\reg^{i+p}_{B, \mc{C}}(S)+\N \mc{C}[1-i-p] \subseteq \mc{Z}^{i+p}
_B (S),
\]
so we can take $\mc{V}_B^{i+p}(S) = \reg^{i+p}_{B, \mc{C}}(S)+\N
\mc{C}[1-i-p]$ in the definition of $ \mc{D}^i _B (J)$, which
shows that this is a special case of Lemma \ref{L:mreg}.
\end{proof}

\begin{Def}
\label{D:reg} Let $J=\{ J_0, J_1, \dots, J_s \}$ with $J_i
\subseteq G$ any subset. For each $i\ge0$, define the
$\N\mc{C}$-modules
\begin{multline*}
\reg^i _{B, \mc{C}}(J)= \{ \mbf{m} \in G:  \mbf{m} + \N \mc{C}[1-i]\\
\subseteq \bigcap_{p = 0}^s \bigcap_{\mbf{d} \in J_p} \mbf{d}
+\reg_{B, \mc{C}}(S) + \N \mc{C}[1-i-p]\}
\end{multline*}
and
\[
\reg _{B, \mc{C}}(J) =  \bigcap_{i \ge 0} \reg^i _{B, \mc{C}}(J).
\]
\end{Def}

\begin{lem}
\label{C:regJ} Suppose $M$ is a $G$-graded $S$-module with a
resolution of type $J$. Then for each $i\geq 0$
\[  \reg^i _{B, \mc{C}}(J) \subseteq
\reg^i_{B, \mc{C}}(M),  \text{ and }  \reg _{B, \mc{C}}(J)
\subseteq \reg_{B, \mc{C}}(M).
\]
\end{lem}

\begin{proof}
If $\mbf {m} \in  \reg^i _{B, \mc{C}}(J) $ then
\[\mbf{m} + \N \mc{C}[1-i] \subseteq
\mbf{d} + \reg_{B, \mc{C}}(S)+\N \mc{C}[1-i-p].
\]
for every $p\ge 0$ and $\mbf{d}\in J_p$. But $\reg_{B,
\mc{C}}(S)\subset \reg^{i+p}_{B, \mc{C}}(S)$, so $\mbf{m}$
satisfies the hypotheses of Lemma \ref{L:regJ}, and the results
follow.
\end{proof}

\begin{lem}
\label{L:nc} Let  $\mc{C} =\{\mbf{c}_1, \dots, \mbf{c}_l\}$.
\begin{itemize}
\item[1.] If $k \le 0$, then
\[
\N \mc{C} [k-1] = \bigcup _{j=1}^{l}(-\mbf{c}_j + \N \mc{C} [k])
\]
\item[2.] If $k \ge 1$, then
\[
\N \mc{C} [k-1] \subseteq \bigcap _{j=1}^{l}(-\mbf{c}_j + \N
\mc{C} [k])
\]
\end{itemize}
\end{lem}

\begin{proof}
Both of these are rather clear.
\end{proof}

 We have a simpler description of $\reg _{B, \mc{C}}(J)$:
\begin{thm}
\label{P:rsimple} Let  $\mc{C} =\{\mbf{c}_1, \dots, \mbf{c}_l\}$.
\[
\reg_{B, \mc{C}}(J) = X \cap Y
\]
where
\begin{align*}
X &=\bigcap_{\mbf{d}\in J_0} ( \mbf{d} + \reg_{B,\mc{C}}(S) +
\N\mc{C}[0])\\
Y&=\bigcap_{p = 1}^s \bigcap_{\mbf{d}\in J_p} \bigcap _{j =
1}^{l}( \mbf{d} - \mbf{c}_j+ \reg_{B, \mc{C}}(S) + \N \mc{C}[1-p])
\end{align*}
\end{thm}

\begin{proof}
We will prove:
\begin{itemize}
\item[1.)] $ \reg_{B, \mc{C}}(J) = \reg_{B, \mc{C}}^0 (J) \cap
\reg_{B, \mc{C}}^1(J) .$ \item[2.)] $ \reg_{B, \mc{C}}^0(J) =
\bigcap_{j=1}^l \bigcap_{p = 0}^s \bigcap_{\mbf{d}\in J_p} (
\mbf{d} - \mbf{c}_j+ \reg_{B, \mc{C}}(S) + \N \mc{C}[1-p]). $
\item[3.)]$ \reg_{B, \mc{C}}^1(J) = \bigcap_{p = 0}^s
\bigcap_{\mbf{d}\in J_p} ( \mbf{d} + \reg_{B,\mc{C}}(S) + \N
\mc{C}[-p]) . $
\end{itemize}
Assuming these, we have
\[ \reg_{B, \mc{C}}(J) = \bigcap_{p = 0}^s \bigcap_{\mbf{d}\in J_p}
X_{p, \mbf{d}} \cap Y_{p, \mbf{d}}, \] where
\begin{align*}
X_{p, \mbf{d}}&= \mbf{d} + \reg_{B,\mc{C}}(S) +\N\mc{C}[-p],\\
Y_{p, \mbf{d}}&= \bigcap_{j=1}^l(\mbf{d} - \mbf{c}_j+ \reg_{B,
\mc{C}}(S) + \N\mc{C}[1-p]).
\end{align*}
But we will see in a moment that $X_{0, \mbf{d}}\subseteq Y_{0,
\mbf{d}}$ and $Y_{p, \mbf{d}}\subseteq X_{p, \mbf{d}}$ if $p \ge
1$. Since $X = \bigcap _{\mbf{d}\in J_{0}}X_{0, \mbf{d}}$ and $Y =
\bigcap _{p= 1}^s\bigcap _{\mbf{d}\in J_{p}}Y_{p, \mbf{d}}$, this
proves the proposition.

The inclusion $X_{0, \mbf{d}}\subseteq Y_{0, \mbf{d}}$ follows
easily from Lemma \ref{L:nc} applied to $k = 1$.

Similarly $Y_{p, \mbf{d}}\subseteq X_{p, \mbf{d}}$ for  $p \ge 1$
follows from the same lemma applied to $k = 1-p \le 0$.

Now we will prove the three claims:

1.) To show $\reg_{B, \mc{C}}(J) = \reg_{B, \mc{C}}^0 (J) \cap
\reg_{B, \mc{C}}^1(J)$, it will suffice to prove that $\reg_{B,
\mc{C}} ^i(J) \subseteq \reg_{B, \mc{C}} ^{i+1}(J)$ for all $i \ge
1$.
\begin{align*}
\mbf{m} \in \reg_{B, \mc{C}} ^i \Leftrightarrow &\ \mbf{m}+\N
\mc{C}[1-i]
\subseteq \mbf{d} + \reg_{B, \mc{C}}(S) +\N \mc{C}[1-i-p], \\
&\forall p\geq 0, \  \ \forall \mbf{d} \in J_p.\\
 \Rightarrow &\  \mbf{m}+\N \mc{C}[1-i] -\mbf{c}_j
\subseteq \mbf{d} + \reg_{B, \mc{C}}(S) +\N \mc{C}[1-i-p]-\mbf{c}_j,\\
&\forall p\geq 0, \  \ \forall \mbf{d} \in J_p,
\ \ \forall j = 1, ..., l.\\
 \Rightarrow &  \  \mbf{m}+\N \mc{C}[1-i] -\mbf{c}_j
\subseteq \mbf{d} + \reg_{B, \mc{C}}(S) +\N \mc{C}[1-(i+1)-p],\\
& \forall p\geq 0, \  \ \forall \mbf{d} \in J_p, \ \ \forall j =
1, ..., l.
\end{align*}
The last implication holds because of Lemma \ref{L:nc}, using that
$i+p \ge 1$. Since $i \ge 1$ we can apply that lemma again to
conclude:
\begin{align*}
 \mbf{m}+\N \mc{C}[1-(i+1)] &= \mbf{m}+
\bigcup _{j = 1}^l(\N \mc{C}[1-i]-\mbf{c}_j)\\
 &=
\bigcup _{j = 1}^l(\mbf{m} + \N \mc{C}[1-i]-\mbf{c}_j )\\
(\mr{from\ the \ above})
 &
\subseteq \mbf{d} + \reg_{B, \mc{C}}(S) +\N \mc{C}[1-(i+1)-p],\\
& \forall p\geq 0, \  \ \forall \mbf{d} \in J_p,\\
\Rightarrow & \ \mbf{m} \in \reg_{B, \mc{C}} ^{i +1}(J).
\end{align*}
2.) We will show $ \reg_{B, \mc{C}}^0(J) = \bigcap_{j=1}^l
\bigcap_{p = 0}^s \bigcap_{\mbf{d}\in J_p} ( \mbf{d} - \mbf{c}_j+
\reg_{B, \mc{C}}(S) + \N \mc{C}[1-p]). $
\begin{align*}
\mbf{m} \in \reg_{B, \mc{C}} ^0(J) \Leftrightarrow &\ \mbf{m}+\N
\mc{C}[1]
\subseteq \mbf{d} + \reg_{B, \mc{C}}(S) +\N \mc{C}[1-p], \\
&\forall p\geq 0, \  \ \forall \mbf{d} \in J_p.\\
 \Leftrightarrow &\  \mbf{m}+\N \mc{C} +\mbf{c}_j
\subseteq \mbf{d} + \reg_{B, \mc{C}}(S) +\N \mc{C}[1-p]\\
&\forall p\geq 0, \  \ \forall \mbf{d} \in J_p,
\ \ \forall j = 1, ..., l.\\
 \Leftrightarrow &  \  \mbf{m}+ \mbf{c}_j
\in \mbf{d} + \reg_{B, \mc{C}}(S) +\N \mc{C}[1-p],\\
& \forall p\geq 0, \  \ \forall \mbf{d} \in J_p, \ \ \forall j =
1, ..., l.
\end{align*}
This last equivalence follows because $\mbf{0}\in \N \mc{C}$ and
each $\N \mc{C}[j]$ is an $\N \mc{C}$-module. This proves the
claim.

3.)Finally, we will prove $ \reg_{B, \mc{C}}^1(J) = \bigcap_{p =
0}^s \bigcap_{\mbf{d}\in J_p} ( \mbf{d} + \reg_{B,\mc{C}}(S) + \N
\mc{C}[-p]) . $
\begin{align*}
\mbf{m} \in \reg_{B, \mc{C}} ^1(J) \Leftrightarrow &\ \mbf{m}+\N
\mc{C}
\subseteq \mbf{d} + \reg_{B, \mc{C}}(S) +\N \mc{C}[-p], \\
&\forall p\geq 0, \  \ \forall \mbf{d} \in J_p.\\
 \Leftrightarrow &  \  \mbf{m}
\in \mbf{d} + \reg_{B, \mc{C}}(S) +\N \mc{C}[-p],\\
& \forall p\geq 0, \  \ \forall \mbf{d} \in J_p, \ \ \forall j =
1, ..., l.
\end{align*}
\end{proof}

Suppose that a region $\mc{D} \subseteq G $ is given. We would
like to define regions $K_p$ such that if $M $ has a resolution of
type $J$ with $J\subseteq K$, that is $J_i \subseteq K_i$ for all
$i$, then $\mc{D} \subseteq \reg_{B, \mc{C}}(M)$.
\begin{Def} \label{D:dreg} For $\mc{D}
\subseteq G$, we define $\dreg_{B,\mc{C}}(\mc{D})= \{ K_0, K_1,
\dots\},$ where
\[ K_p=\{ \mbf{d}\in G: \mc{D} +\N\mc{C}[1-i]
\subseteq  \mbf{d}+\reg_{B, \mc{C}}(S) + \N\mc{C}[1-i-p], \
 \ \forall i\geq 0\}. \]  We also denote
\[ K_{p,i}=\{ \mbf{d}\in G: \mc{D} +\N\mc{C}[1-i]
\subseteq  \mbf{d}+\reg_{B, \mc{C}}(S) + \N\mc{C}[1-i-p]\},\]
  and
 $\dreg_{B,\mc{C}}(\mc{D})_p=K_p=\bigcap_{i\geq 0} K_{p,i}$.
\end{Def}

We drop $B$ or $\mc{C}$ from the notation if it is clear.
\begin{lem}
\label{L:dreg} Suppose that  $\mc{D} \subseteq G $ is given and
that $M$ is a $G$-graded $S$-module of type $J$ with $J \subseteq
\dreg_{B, \mc{C}}(\mc{D}) $. Then $\mc{D} \subseteq \reg_{B,
\mc{C}}(M)$.
\end{lem}

\begin{proof}
By definition, if $\dreg_{B,\mc{C}}(\mc{D})= \{ K_0, K_1,
\dots\}$, then
\[ \mc{D}+\N \mc{C}[1-i] \subseteq \bigcap_{\mbf{d} \in K_p}
( \mbf{d} + \reg(S) + \N \mc{C}[1-i-p]), \  \ \forall p\geq 0.\]
Since $J \subseteq K$, by which we mean $J_p \subseteq K_p$ for
all $p$,
\[ \bigcap_{\mbf{d} \in K_p}
( \mbf{d} + \reg(S) + \N \mc{C}[1-i-p])
 \subseteq \bigcap_{\mbf{d} \in J_p}
( \mbf{d} + \reg(S) + \N \mc{C}[1-i-p]).
\]
Combining these with Lemma \ref{L:regJ} gives the result.
\end{proof}

\begin{lem} \label{L:dmod} Let $\bar{\mc{D}}$ be the $\N \mc{C}$-module
generated by $\mc{D}$. Then
\[
\dreg_{B, \mc{C}}(\mc{D})=\dreg_{B, \mc{C}}(\bar{\mc{D}}).
\]
\end{lem}

\begin{proof}
It is true since $\N \mc{C}[i]$ is an $\N \mc{C}$-module
 for each $i$.
\end{proof}

We can give a simpler description of these regions:

\begin{thm}
\label{P:drsimple} Let  $\mc{C} =\{\mbf{c}_1, \dots, \mbf{c}_l\}$.
Then $\mbf{d} \in \dreg_{B, \mc{C}}(\mc{D})_p$ if and only if
\[
\mc{D} \subseteq (\mbf{d} + \reg _{B, \mc{C}}(S) +\N \mc{C}[-p] )
\cap \bigcap _{j = 1}^{l} (\mbf{d} - \mbf{c}_j + \reg _{B,
\mc{C}}(S) +\N \mc{C}[1-p]).
\]
\end{thm}

\begin{proof}
By Definition \ref{D:dreg}, $\dreg_{B,\mc{C}} (\mc{D})_p = K_p =
\bigcap _{i \ge 0}K_{p, i}.$   If we can show $K_{p, i} \subseteq
K_{p, i+1}$ for $i\ge 1$, then $K_p=K_{p, 0} \cap K_{p, 1}$.
Namely,
\[
\mbf{d} \in K_{p, i} \Leftrightarrow \mc{D} + \N \mc{C}[1-i]
\subseteq \mbf{d} + \reg _{B, \mc{C}}(S) +\N \mc{C}[1-i-p].
\]
When $i = 0$ this is equivalent to
\[
\mbf{d} \in K_{p, 0} \Leftrightarrow \mc{D} + \mbf{c}_j \subseteq
\mbf{d} + \reg _{B, \mc{C}}(S) +\N \mc{C}[1-p].
\]
for all $j = 1, ..., l$, since each $\N \mc{C}[k]$ is an $\N
\mc{C}$-module. Or:
\[
\mbf{d} \in K_{p, 0} \Leftrightarrow \mc{D} \subseteq \bigcap _{j
= 1}^{l} (\mbf{d} - \mbf{c}_j + \reg _{B, \mc{C}}(S) +\N
\mc{C}[1-p]).
\]
Similarly, when $i=1$,
\[
\mbf{d} \in K_{p, 1} \Leftrightarrow \mc{D} \subseteq \mbf{d} +
\reg _{B, \mc{C}}(S) +\N \mc{C}[-p].\]
 Now we show $K_{p,i}
\subseteq K_{p, i+1}$ for $i\geq 1$.
\begin{align*}
\mbf{d} \in K_{p,i} \Leftrightarrow &\ \mc{D}+\N \mc{C}[1-i]
\subseteq \mbf{d} + \reg_{B, \mc{C}}(S) +\N \mc{C}[1-i-p]\\
 \Rightarrow &\  \mc{D}+\N \mc{C}[1-i] -\mbf{c}_j
\subseteq \mbf{d} + \reg_{B, \mc{C}}(S) +\N \mc{C}[1-i-p]-\mbf{c}_j,\\
& \forall j = 1, ..., l.\\
 \Rightarrow &  \  \mc{D}+\N \mc{C}[1-i] -\mbf{c}_j
\subseteq \mbf{d} + \reg_{B, \mc{C}}(S) +\N \mc{C}[1-(i+1)-p],\\
& \forall j = 1, ..., l.
\end{align*}
The last implication holding because of Lemma \ref{L:nc} since
$1-i-p \le 0$. Since $i \ge 1$ we can apply that lemma again to
conclude:
\begin{align*}
 \mc{D}+\N \mc{C}[1-(i+1)] &= \mc{D}+
\bigcup _{j = 1}^l(\N \mc{C}[1-i]-\mbf{c}_j )\\
 &\subseteq
\bigcup _{j = 1}^l(\mc{D} + \N \mc{C}[1-i]-\mbf{c}_j )\\
(\mr{from\ the \ above})
 &
\subseteq \mbf{d} + \reg_{B, \mc{C}}(S) +\N \mc{C}[1-(i+1)-p].
\end{align*}
By Definition \ref{D:dreg}, this implies that $\mbf{d} \in K_{p, i
+1}$.  Therefore, we have $\mbf{d} \in \dreg_{B,
\mc{C}}(\mc{D})_p=K_p$ if and only if
\[
\mc{D} \subseteq (\mbf{d} + \reg _{B, \mc{C}}(S) +\N \mc{C}[-p] )
\cap \bigcap _{j = 1}^{l} (\mbf{d} - \mbf{c}_j + \reg _{B,
\mc{C}}(S) +\N \mc{C}[1-p]).
\]
\end{proof}

The above results imply that if the region $\mc{D}$ is finitely
generated as an  $\N \mc{C}$-module the conditions defining
$\dreg(\mc{D})$ are finite in number - a priori, there are only
finitely many possible sets $K_p$ to be considered in the free
resolutions of any module.

The following is easily verified:
\begin{cor}
\label{L:corresp} The mappings
\[ \left (\begin{matrix}
\mr{collections}\\
K = \{ K_0, K_1, \dots, K_s\} \\
K_i \subseteq G \end{matrix} \right)
\begin{matrix} \reg \\ \rightarrow \\ \leftarrow \\ \dreg
\end{matrix}
\left( \begin{matrix}
\N \mc{C}-\mr{modules}\\
\mc{D} \subseteq G \\  \end{matrix} \right)\] satisfy
\begin{itemize}
\item[1.] If $K\subseteq K'$ denotes $K_i \subseteq K'_i$ for all
$i\geq 0$, then $\reg_{B, \mc{C}}(K) \supseteq \reg_{B,
\mc{C}}(K').$ \item[2.] If $\mc{D} \subseteq \mc{D'}$, then
$\dreg_{B, \mc{C}}(\mc{D}) \supseteq \dreg_{B, \mc{C}}(\mc{D'})$.
\item[3.] $K \subseteq \dreg_{B, \mc{C}} (\reg_{B, \mc{C}}(K))$.
\item[4.] $\mc{D} \subseteq \reg_{B, \mc{C}}(\dreg_{B,
\mc{C}}(\mc{D}))$.
\end{itemize}
\end{cor}

\section{\textbf{v}-gradings}\label{S:vdot}
Let $S=K[x_1, \dots, x_n]$ be a $G$-graded algebra where $G=\Z^r$.
\begin{Def}
\label{D:vgrade} A vector $\mbf{v} \in \Z^r$ is called a {\it
coarsening vector} for $S$ if
\[\deg_{\mbf{v}}(x_i) =\mbf{v} \cdot \deg(x_i) \geq 0,  \  \
\forall 1\leq i \leq n.\] It is called a {\it positive coarsening
vector} if $\deg_{\mbf{v}}(x_i) >0$ for all $i$.
\end{Def}
The ring $S$ and all $G$-graded $S$-modules $M$ inherit
$\Z$-gradings from $\mbf{v}$.  We denote these $\Z$-graded objects
by $S_{\mbf{v}}$ and $M_{\mbf{v}}$, where
\[ (M_{\mbf{v}})_{m} = \bigoplus_{\mbf{d}\in G, \mbf{v} \cdot
\mbf{d} = m } M_{\mbf{d}}.\] We can regard the variables $x_i$
with $\mbf{v} \cdot \deg(x_i)=0$ as constants.  That is we can
write $S_{\mbf{v}}=R[x_j]$ where $x_j$ with $\mbf{v} \cdot
\deg(x_i) > 0$, and $R=K[x_i]$ with $\mbf{v} \cdot \deg(x_i) = 0$.
Then $S_{\mbf{v}}$ is the homogeneous coordinate ring of the
weighted projective space over $R$ in the coordinates $x_j$.  If
the vector $\mbf{v}$ is understood, we will simply write $S$ and
$M$, and the notation $M_m$ for $m\in \Z$ will refer to this
$\mbf{v}$-grading.

Assume $\mbf{v}$ is a positive coarsening vector. Let
$c_{\mbf{v}}=\lcm \{\deg_{\mbf{v}}(x_i): 1\leq i \leq n\}.$ Recall
that $\N c_{\mbf{v}}$ is the Kaehler cone for $S_{\mbf{v}}$, for
details see \cite[Example 2.1]{MS}.

\begin{Def}
\label{D:vreg} If $M$ is a finitely generated $G$-graded
$S$-module and $\mbf{v}$ is a positive coarsening vector for $S$,
then
\begin{equation*}
\reg_{\mbf{v}}(M) = \{ p \in \Z: H^i_{\mfr{m}}(M)_q = 0, \  \
\forall q \in p+ \N c_{\mbf{v}}[1-i], \mfr{m}=(x_1, \dots, x_n)\}.
\end{equation*}
This is $\mfr{m}$-regularity with respect to $\mc{C}=\{
c_{\mbf{v}}\}$ defined by \cite{MS} for the $\mbf{v}$-graded
module $M$.
\end{Def}

Here we will quote Theorem \ref{T:numv} and Theorem \ref{T:resnum}
from \cite{SVW}.
\begin{thm}\label{T:numv}
Let $M \neq 0$ be a finitely generated $G$-graded module. Let
${\mbf v} \in G$ be a positive coarsening vector for $S$. Then
there exists a $p \in \reg_{\mbf{v}}(M)$ such that $q \ge p$
implies that $q \in \reg_{\mbf{v}}(M)$. The least such $p$ with
this property is called $\vregnum(M)$. If $M = 0$, then we set
$\vregnum(M) = -\infty$.
\end{thm}

\begin{thm}\label{T:resnum}
Let $S,  M, \mbf{v}$ be as in the previous theorem.  Let
$C_{\bullet}$ be a minimal $G$-graded free resolution of $M.$ Let
\[s_{\mbf{v}} = \max\{nc_{\mbf{v}} - \sum \deg_{\mbf{v}}(x_i), c_{\mbf{v}}\}.\]
Then the degrees $\mbf{d} \in G$ of the module $C_i$ satisfy
\[
\deg_{\mbf{v}}(\mbf{d}) \leq
\vregnum(M)+is_{\mbf{v}}+c_{\mbf{v}}-1.
\]
\end{thm}

\begin{Def}
\label{Dbstarreg} Let $\mbf{v}_1, \dots, \mbf{v}_t$ be positive
coarsening vectors for $S$.  Let $b_1, \dots, b_t$ be integers.
Let $M$ be a finitely generated $G$-graded $S$-module.  We say
that $M$ is {\it $b_*$-regular with respect to $\mbf{v}_*$} if
\[ b_j \in \reg_{\mbf{v}_j}(M), \ \ \forall 1 \leq j \leq t.\]
\end{Def}

\begin{cor}
\label{C:resnum} Suppose that $M$ is $b_*$-regular with respect to
$\mbf{v}_*$, where each $b_j \geq \vjregnum(M)$.  Then the degrees
$\mbf{d}$ of the $i$-th syzygy module $C_i$ in a minimal
$G$-graded resolution of $M$ belong to the convex polyhedral
region
\[K_i(\mbf{v}_*, b_*) =\{ \mbf{d} \in G: \deg_{\mbf{v}_j}(\mbf{d})
\leq b_j+is_{\mbf{v}_j}+c_{\mbf{v}_j} - 1,  \  \ \forall 1\leq j
\leq t\}.\]
\end{cor}

\begin{rem}\label{R:ki}
Let $Q \subseteq G$ be the sub-semigroup generated by the degrees
$Q=\N\{\deg(x_i): 1\leq i \leq n\}$.  Then the support of $M$ is
contained in a union of finitely many translates $\mbf{b}_k+Q$.
The supports of the syzygy module $C_i$ in a minimal $G$-graded
resolution of $M$ will be contained in the same union of these
translates, so we can actually say that under the hypotheses of
Theorem \ref{T:resnum}, the degrees of $C_i$ are in the finite set
\[ (\bigcup_{j=1}^k(\mbf{b}_j+Q) )\bigcap K_i(b_*, \mbf{v}_*),
\  \ \text{for some } k.\]
\end{rem}

\begin{thm}
\label{P:vectres} Suppose that $M$ is a finitely generated
$G$-graded $S$-module, $\mbf{v}_1, \dots, \mbf{v}_t \in \Z^r$
positive coarsening vectors for $S$, and $b_1, \dots, b_t \in \Z$
integers such that $b_j \geq \vjregnum(M)$ for  $1 \leq j \leq t$.
Suppose that $M$ is $b_*$-regular with respect to $\mbf{v}_*$. Let
$J = \{ J_0, J_1, ..., J_s\}$ where
\[
J_p = (\bigcup_{j=1}^k( \mbf{b}_j+  \N \mc{C} ) ) \bigcap
K_p(\mbf{v}_*,b_*).
\]
Then,
\begin{itemize}  \item[1.] $M$ has a resolution of
type $ J'$,  with $J'=\{J'_0, \dots, J'_s\}$ and $J'_i\subseteq
J_i$. \item[2.] $\reg(J) \subseteq \reg(J') \subseteq \reg_{B,
\mc{C}} (M)$.
\end{itemize}
\end{thm}

\begin{proof}
The first claim is a restatement of Remark \ref{R:ki}. The second
claim follows from the first claim and by Lemma \ref{C:regJ}.
\end{proof}

The above gives a lower bound on a regularity region for $M$ in
terms of a $b_*$-regularity region with respect to a collection of
coarsening vectors. We now give a set of conditions, based on
vanishing of $\lclcoh{i}{B}{M}_{\mbf{d}}$ for $\mbf{d}$ in half
planes, that imply $\mbf{m}\in \reg_{B, \mc{C}} (M)$.

\begin{Def}
\label{D:halfplane} Let $\mbf{v} \in \R^r$, $b\in \R$, define the
half-planes and the hyperplane
\begin{align*}
& P^+_{\mbf{v}, b}=\{ \mbf{x} \in \R^r : \mbf{v} \cdot \mbf{x}
\geq b\}, \  \ P^-_{\mbf{v}, b}=-P^+_{\mbf{v}, b}, \\
&L_{\mbf{v}, b}= P^-_{\mbf{v}, b}\cap P^+_{\mbf{v}, b}=\{\mbf{x}
\in \R^r : \mbf{v} \cdot \mbf{x} = b\}.
\end{align*}
\end{Def}

\begin{lem}\label{L:PL}
Let $\mc{C} =\{\mbf{c}_1, \dots, \mbf{c}_l\}$. Let $\mbf{v} \in
\R^r$.  Suppose that
\[
\N \mc{C} \subseteq P^+_{\mbf{v}, 0}, \ \text{that\  is \ }
\mbf{v} \cdot \mbf{c}_j \geq 0 \text{\ for\ all\ } j=1, \dots, l.
\]
Note that such a vector $\mbf{v}$ always exists under the
assumption $\mc{C}\subset\mc{K}^{\mr{sat}} $ since the latter is a
subsemigroup of $\N \mc{A}$ and we are assuming that
$\mr{pos}(\overline{A})$ is a pointed cone (see the introduction).
If we define $ \min(\mbf{v}, \mbf{c}) = \min(\mbf{v} \cdot
\mbf{c}_j)_{j=1}^l$, and $\max(\mbf{v}, \mbf{c}) = \max(\mbf{v}
\cdot \mbf{c}_j)_{j=1}^l.$ Then $\N \mc{C}[k] \subseteq
P^+_{\mbf{v}, b}$, where $b=k\min(\mbf{v}, \mbf{c})$ if $k \geq
0$, and $b= k \max(\mbf{v}, \mbf{c})$ if $k <0$.
\end{lem}

\begin{proof}
For $k=0$ it is true by assumption. Observe that
\begin{align*}
\N \mc{C}[k+ 1]= \bigcup_{j=1}^l
\mbf{c}_j+ \N \mc{C}[k], \text{ if  } k \ge 0,\\
\N \mc{C}[k- 1]= \bigcup_{j=1}^l -\mbf{c}_j+ \N \mc{C}[k], \text{
if  } k\le 0.
\end{align*}
We will prove the result by performing induction on $k$.  Assume
that $k \ge 0$. If $\mbf{d} \in\N \mc{C}[k+ 1] $, then $\mbf {d} =
\mbf{c}_j + \mbf{e}$ for some $\mbf{e} \in  \N \mc{C}[k]$. By
induction hypothesis, we have \[\mbf{v}\cdot\mbf {d} = \mbf{v}
\cdot (\mbf{c}_j + \mbf{e}) \ge \mbf{v}\cdot \mbf{c}_j + k \min
(\mbf{v}, \mbf{c}) \ge
 (k+1) \min (\mbf{v}, \mbf{c}).\]  By Definition \ref{D:halfplane},
 we have that
$\N \mc{C}[k] \subseteq P^+_{\mbf{v}, b}$ where $b=k\min(\mbf{v},
\mbf{c})$.

Now, assume that $k \le 0$. If $\mbf{d} \in\N \mc{C}[k-1] $, then
$\mbf {d} = -\mbf{c}_j + \mbf{e}$ for some $\mbf{e} \in  \N
\mc{C}[k]$. By induction hypothesis, we have \[\mbf{v}\cdot\mbf
{d} = \mbf{v} \cdot (-\mbf{c}_j + \mbf{e} ) \ge \mbf{v}\cdot
(-\mbf{c}_j) + k \max (\mbf{v}, \mbf{c}) \ge
 (k-1) \max (\mbf{v}, \mbf{c}).\]  By Definition \ref{D:halfplane},
 we have that
$\N \mc{C}[k] \subseteq P^+_{\mbf{v}, b}$ where $b=k\max(\mbf{v},
\mbf{c})$.
\end{proof}

\begin{thm}
\label{L:degbound} Let $\mbf{v} \in G$ be a coarsening vector for
$S$ such that $\mbf{v} \cdot \mbf{c}_j \geq 0$ for all $j =1,
\dots, l$.  Let $k\in \Z$. Let $M$ be a $G$-graded $S$-module.
Fix $\mbf{m}\in G$, and let $B\subseteq G$ be a $G$-graded ideal.
Suppose that
\[H^i_B(M)_{\mbf{d}}=0, \   \ \forall \mbf{d}, \text{ such that }
\deg_{\mbf{v}}(\mbf{d}) \geq\deg_{\mbf{v}}(\mbf{m})  + (k-i) \cdot
m(k, i),\] where \[ m(k,i)= \max(\mbf{v}, \mbf{c}), \text{ if }
k-i < 0, \text{ and } m(k,i)= \min(\mbf{v}, \mbf{c}), \text{ if }
k-i \geq 0.\] Then
\[\mbf{m} +\N \mc{C}[k-i] \subseteq
\mc{Z}^i_B(M).\] In particular, if $k=1$ and this holds for all $i
\geq 0$, this says that $\mbf{m} \in \reg_{B, \mc{C}}(M)$.
\end{thm}

\begin{proof}
The condition states that
\[H^i_B(M)_{\mbf{d}}=0, \   \
\forall \mbf{d}-\mbf{m} \in P^+_{\mbf{v}, b}, \text{ where }
b=(k-i) \cdot m(k,i). \] By Lemma \ref{L:PL}, $\N \mc{C}[k-i]
\subseteq P^+_{\mbf{v}, b}$. Therefore
\[H^i_B(M)_{\mbf{d}}=0, \   \
\forall \mbf{d}-\mbf{m} \in \N\mc{C}[k-i],  \] that is
\[\mbf{m} +\N \mc{C}[k-i] \subseteq
\mc{Z}^i_{B, \mc{C}}(M).\] By 2 and 3 in Definition \ref{D:reg},
when $k=1$ and all $i \geq 0$, we have $\mbf{m} \in \reg_{B,
\mc{C}}(M)$.
\end{proof}

\section{Regularity with respect to a family of ideals}\label{S:family}

Let $(S, B)$ be a $G$-graded $R$-algebra with a monomial ideal. In
this section we consider the case where there is a decomposition
\begin{align*}
& Z=\V(B)=\bigcup_{i=1}^t \V(B_i)= \bigcup_{i=1}^t Z_i\\
&W_I = \bigcap_{i \in I} Z_i =\V(B_I), \  \ \emptyset \neq I
\subseteq [t] = \{1, \dots, t\}
\end{align*}
If $(S, B)$ is the homogeneous coordinate ring and irrelevant
ideal of a toric variety $\proj{}_{\Sigma}$ (see \cite{C}), then
there is a canonical decomposition of this type, due to Batryev
(see \cite[section 10]{C2}) in which the ideals $B_I$ are defined
by linear forms. Here is the definition: The cones of the fan
$\Sigma$ are spanned by integral vectors in a vector space
$N_{\R}$. A collection of edge generators $ \{ \mbf{n}_{i_{1}},
..., \mbf{n}_{i_{r}}\} $ is called primitive if they do not all
lie in a cone of $\Sigma$ but every subset of it does. Then, we
have a decomposition as above with
\[
B_i = \langle  x_{i_{1}}, ..., x_{i_{r}}\rangle, \ B_I = \sum
_{i\in I}B_i
\]
where the index $i$ runs over all the primitive sets.

\begin{lem}
\label{L:famvan} For each $i\ge 0$
\[\bigcap_{p\geq 0}
\bigcap_{\substack {I\subseteq [t]\\ \#I=p+1 }}
\mc{Z}_{B_I}^{i+p}(M) \subseteq \mc{Z}_B^i(M).\]
\end{lem}
\begin{proof}
Let $\mbf{d} \in \bigcap_{p\geq 0} \bigcap_{\substack {I\subseteq
[t]\\ \#I=p+1 }} \mc{Z}_{B_I}^{i+p}(M)$, then $\mbf{d} \in G$ and
$H^{i+p}_{B_I}(M)_{\mbf{d}}=0$ for all $p\geq 0$ and $I\subseteq
[t]$ with $\#I=p+1$ by Definition \ref{D:reg}.  By Corollary
\ref{C:1}, this implies that $H^i_{B}(M)_{\mbf{d}}=0$, that is,
$\mbf{d} \in \mc{Z}_B^i(M)$.
\end{proof}
\begin{lem}\label{L:minregb}
Fix $i\ge 0$. If \[ \mbf{m} + \N \mc{C}[1-i] \subseteq
\bigcap_{p\geq 0} \bigcap_{\substack {I\subseteq [t]\\ \#I=p+1 }}
\mc{Z}_{B_I}^{i+p}(M),\] then $\mbf{m} \in \reg^i_{B, \mc{C}}(M)$.
If this is true for all $i$, then
 $\mbf{m} \in \reg_{B, \mc{C}}(M)$.
\end{lem}

\begin{proof}
By Lemma \ref{L:famvan}, this implies
\[ \mbf{m} + \N \mc{C}[1-i] \subseteq   \mc{Z}_B^i(M). \]
By Definition \ref{D:reg}, $\mbf{m} \in\reg^i_{B, \mc{C}}(M) $;
and $\mbf{m} \in \reg_{B, \mc{C}}(M)$ if this is true for all $i$.
\end{proof}

\begin{Def}
\label{D:regstar} A $G$-graded $S$-module $M$ is $\mbf{m}$-regular
in dimension $i$ with respect to the collection $\{B_I\}$ and
$\mc{C}$ if for all
\[ \mbf{m} + \N \mc{C}[1-i] \subseteq  \mc{Z}_{B_I}^{i+\#I-1}(M),
\text{\ for\ all\ } \emptyset \neq I \subseteq [t]. \] We denote
this by $\mbf{m} \in \reg^i_{B_*, \mc{C}}(M)$, and define
$$\reg_{B_*, \mc{C}}(M) =\bigcap_{i\geq 0} \reg^i_{B_*, \mc{C}}(M).$$
\end{Def}

\begin{pro}
\label{P:regstar} Suppose that $B_i$ is a family of monomial
ideals of $S$ giving a decomposition of a monomial ideal $B$ as
indicated above. Then
\[\reg_{B_*, \mc{C}}(M) \subseteq \reg_{B, \mc{C}}(M).\]
\end{pro}

\begin{proof}
If $ \mbf{m} \in \reg_{B_*, \mc{C}}(M)$, then by Definition
\ref{D:regstar},
\[ \mbf{m}+ \N \mc{C}[1-i]\subseteq
\bigcap_{p\geq 0} \bigcap_{\substack {I\subseteq [t]\\ \#I=p+1 }}
 \mc{Z}_{B_I}^{i+p}(M),  \   \ \forall i\geq 0.\] By Lemma
\ref{L:minregb}, we have that
$\mbf{m} \in \reg^i _{B, \mc{C}}(M)$ for all $i \geq 0$.
Therefore, $\mbf{m} \in \bigcap_{i\geq 0} \reg^i _{B, \mc{C}}(M)
=\reg _{B, \mc{C}}(M)$.
\end{proof}

Assume $G= \Z ^r$. Let us consider the case of a decomposition of
$Z = \V (B)$. When discussing vanishing conditions for
$\lclcoh{i}{B_{I}}{M}$ it is natural to consider coarsening
vectors that assign $\mbf{0}$ to some of the variables. This is
justified, that is, that one can consider coarsenings that assign
zero to some of the variables, will be shown in the Section
\ref{S:0degree} (see especially Remark \ref{R:application}). In
fact for each subset $I$ one could consider a finite set of such
coarsenings. We will do the simplest case of assigning only
one coarsening vector to each $\emptyset \ne I$.

\begin{Def}\label{D:Bv}
Let $S$ be a $G=\Z ^r$-graded $R$-algebra. Let $B_i, \ i\in [t]$
be a family of monomial ideals in $S$. For each $\emptyset \ne I
\subseteq [t]$, let $\mbf{v}_{I}\in G$ be a coarsening vector, not
necessarily strict. Let $M$ be a $G$-graded $S$-module. We say
that $\mbf{m}\in \reg ^i _{B_{*}, \mbf{v}_{*}, \mc{C}}$ if for all
$\emptyset \ne I \subseteq [t]$,
\[
\lclcoh{i}{B_{I}}{M}_{\mbf{d}}=0\ \text{for\ all\ } \deg
_{\mbf{v}_{I}} (\mbf{d})\ge\deg _{\mbf{v}_{I}} (\mbf{m})
                 +(1 - i)m(1 , i)
\]
where
\[ m(1,i)= \max(\mbf{v}_I, \mbf{c}), \text{ if } 1-i <
0, \text{ and } m(1,i)= \min(\mbf{v}_I, \mbf{c}), \text{ if } 1-i
\geq 0.\] Define
\[
\reg_{B_*,\mbf{v}_*, \mc{C}}(M) =\bigcap_{i\geq 0}
\reg^i_{B_*,\mbf{v}_* ,\mc{C}}(M).
\]
\end{Def}
\begin{pro}
\label{P:Bv} Notations as in Definition \ref{D:Bv}, if
$\mbf{v}_I\cdot \mbf{c}_j \ge 0$ for all $\emptyset \ne I
\subseteq [t]$ and all $ \mbf{c}_j\in \mc{C}$, then
\[\reg_{B_*,\mbf{v}_*, \mc{C}}(M) \subseteq
\reg_{B_*, \mc{C}}(M).\] If the ideals $B_i$ arise from a
decomposition of $\V (B)$ as explained at the beginning of this
section, then
\[\reg_{B_*,\mbf{v}_*, \mc{C}}(M)
\subseteq\reg_{B,\mc{C}}(M).\]
\end{pro}
\begin{proof}
If $\mbf{m}\in \reg_{B_*,\mbf{v}_*, \mc{C}}(M) $ then for all
$\emptyset \ne I \subseteq [t]$, the conditions of Theorem
\ref{L:degbound} hold with $B = B_I$, $\mbf{v} = \mbf{v}_I$ and $k
= 1$. From the result of that Theorem, we conclude that
\[
\mbf{m} +\N \mc{C}[1-i] \subseteq \mc{Z}^{i}_{B_I}(M), \ \text{for
\ all\ } \ i \ge 0.
\]
In particular,
\[
\mbf{m} +\N \mc{C}[(2-\# I) -i] \subseteq \mc{Z}^{i+\# I -
1}_{B_I}(M), \ \text{for \ all\ } \ i \ge 0.
\]
But $(2 - \# I) - i\le 1 -i $, thus
\[
\mbf{m} +\N \mc{C}[1 -i] \subseteq\mbf{m} +\N \mc{C}[(2-\# I) -i]
\subseteq  \mc{Z}^{i+\# I - 1}_{B_I}(M), \text{ for all } i\geq 0.
\]
By Definition \ref{D:regstar}, $\mbf{m} \in \bigcap_{i\geq 0}
\reg^i_{B_*,\mc{C}}(M)= \reg_{B_*,\mc{C}}(M)$.

If we assume that the $B_*=\{B_i\}$ comes from a decomposition of
a monomial ideal $B$, then $\mbf{m} \in \reg_{B_*,\mc{C}}(M)
\subseteq \reg_{B,\mc{C}}(M)$ by Proposition \ref{P:regstar}.
\end{proof}

Often it is the case that if $\mbf{m}\in \reg_{B_*, \mbf{v}_* ,
\mc{C}}(M)$, we can find a resolution for $M$ with syzygies in a
particular region depending on $\mbf{m}$. The point of considering
a family of ideals is that imposing regularity conditions with
respect to each ideal $B_I$ gives a family of constraints, indexed
by $I$, on the degrees of syzygies of a minimal resolution of a
module. To apply the methods of Section \ref{S:vdot} here, it is
necessary to assume that the coarsenings $\mbf{v}_I$ are
orthogonal in the following sense:

Assume that each $B_I$ is generated by a subset of generators
$x_j$ of $S = K[x_1, ..., x_n]$. Then we must postulate that $\deg
_{\mbf{v}_I}(x_j) > 0$ for every variable $x_j \in B_I$, and that
$\deg _{\mbf{v}_I}(x_j) =0$ for every variable $x_j \notin B_I$.
Such a vector may not exist, in general, but they do exist for
fans of the form $\Sigma = \Sigma _1 \times ... \times \Sigma _r$.

For instance, if we assume that the variables in $B_i$ are
orthogonal to those of $B_j$ for $i \ne j$ in the sense that we
can choose coarsenings $\mbf{v}_i$ for the variables $x_k \in B_i$
with $\deg _{\mbf{v}_i}(x_k) > 0$ for every variable $x_k \in
B_i$, and  $\deg _{\mbf{v}_I}(x_k) =0$ for every variable $x_k
\notin B_j$ , we can define $\mbf{v}_I =\sum_{i \in I}\mbf{v}_i$
since $B_I = \sum _{i \in I}B_i$.

\begin{pro}
\label{P:vres} Let $B_i, \ i\in [t]$ be the set of ideals arising
from the Batryev decomposition of a toric variety, as in the
introduction of this section. Recall that each ideal $B_I$ is
generated by a subset of generators $\{x_1, ..., x_n\}$ of $S =
K[x_1, ..., x_n]$. Suppose that coarsening vectors are chosen
$\mbf{v}_I$ for $\emptyset \ne I \subseteq [t]$ with the property
that $\deg _{\mbf{v}_I}(x_j) > 0$ for every variable $x_j \in B_I$
and that  $\deg _{\mbf{v}_I}(x_j) =0$ for every variable $x_j
\notin B_I$.

Let $M$ be a $G = \Z ^r$-graded $S$-module and $\mbf{m} \in
\reg_{B_*,\mbf{v}_* ,\mc{C}}(M)$. Suppose for each $I$ integers
$b_I$ are given with $b_I\ge \vIregnum (M)$.

Assume that for every $I$, and $i = 0, ..., n$,
\[
\deg _{\mbf{v}_I} (\mbf{m})\le b_I+(1-i)(c_{\mbf{v}_I} - m(1, i)).
\]
If we let $\mbf {v}_*$ be the collection $\{\mbf {v}_I \}$ and
$b_*$ be the collection of integers $\{b_I\}$, then $M$ has a
resolution of type $J'$ where $J'$ is defined as in Theorem
\ref{P:vectres}.
\end{pro}
\begin{proof}
The condition that $\mbf{m} \in \reg_{B_*,\mbf{v}_* ,\mc{C}}(M)$
and the given bounds on $\deg _{\mbf{v}_I} (\mbf{m})$ imply that
for all $\mbf{d}$ with
\[
\deg _{\mbf{v}_I}(\mbf{d}) \ge b_I + (1-i)c_{\mbf{v}_I}.
\]
Therefore, we have $\lclcoh{i}{B_{I}}{M}_{\mbf{d}}= 0$. Since
$b_I\ge \vIregnum (M)$, which holds for all $I$, this is
equivalent to the condition that $M$ is $b_*$-regular with respect
to $\mbf {v}_*$.  Thus we may apply Theorem \ref{P:vectres} to get
the result.
\end{proof}

\section{regularity for $G$-graded $R$-algebras}
\label{S:0degree} In this section we show how the theory developed
so far can be extended from a $G$-graded polynomial ring $S =
K[x_1, ..., x_n] $ where $K$ is a field to $S = R[x_1, ..., x_n]$
where $R$ is a commutative Noetherian ring. If $S = R[x_1, ...,
x_n]$, a finitely generated $G$-graded $S$-module may not have a
finite free graded resolution of the type we have been
considering. In our applications this is no restriction: we always
start with a module over $S = K[x_1, ..., x_n] $ which a priori
has such free $G$-graded finite resolution, and we seek bounds on
the degrees of generators of syzygies by consideration of vector
gradings that possibly assign degree zero to some of the variables
(see Remark \ref{R:application}).

Note that all the results of section \ref{S:relation}
are true as stated, with the ground field $K$ replaced
by any commutative ring $R$. First we define regularity
as in the introduction, but for $S= R[x_1, ..., x_n]$.
Then the results of section \ref{S:relation} follow
formally from the spectral sequence arguments of section
\ref{S:tool}: we are postulating the free resolutions
for our modules.

The results of section \ref{S:vdot} depend on the paper
\cite{SVW}, whose main results depend on Maclagan and Smith's
paper. We will argue that the main results,
especially \cite[Theorem 4.7, Theorem 5.4]{MS} are true as stated
in this more general setting.

This is done in two stages:
\begin{itemize}
\item[a.] The theory works for $R$, a local ring with infinite
residue field. \item[b.] If the theory works for all local rings
with infinite residue field, then it works for all commutative
Noetherian rings.
\end{itemize}
By ``theory works'' we mean that the results stated in \cite{MS}
and \cite{SVW} are valid for the ring $S = R[x_1, ..., x_n]$ where
$R$ is the type of ring under discussion. In all cases, $B\subset
S$ is a monomial ideal, possibly subject to further restrictions
as needed in those cited works. The main results of those works
are of the following sort:
\begin{itemize}
\item[1.] Knowing that $\lclcoh{i}{B}{M}_{\mbf{d}}= 0$ for
$\mbf{d}$ in a certain region of $G$, deduce that
 $\lclcoh{i}{B}{M}_{\mbf{d}}= 0$ for
$\mbf{d}$ in a larger region. \item[2.] From the vanishing
$\lclcoh{i}{B}{M}_{\mbf{d}}= 0$ for $\mbf{d}$ in certain regions
of $G$ depending on $i$, deduce that the degrees of the generators
of $S$-module $M$ may be found in a region of $G$. More generally
the degrees of the $p$th syzygies are in a certain region.
\end{itemize}

Let us first show how assuming we have verified step a. above
we may conclude step b. The idea,
due to Ooishi \cite{O}, is to consider the localizations $R \to
R_{\mfr{p}}\to  R_{\mfr{p}}(T)$ where $\mfr{p}$ ranges over $\Spec
(R)$ and $T$ is a variable. The first arrow is flat and the second
arrow is faithfully flat. Thus, if $P\subset Q$ are $R$-modules,
then $P=Q$ if and only if $P_{\mfr{p}}(T) = Q_{\mfr{p}}(T)$ for
all $\mfr{p}\in \Spec (R)$. Applied to $P = 0$, $Q =
\lclcoh{i}{B}{M}_{\mbf{d}}$  we see that
$\lclcoh{i}{B}{M}_{\mbf{d}}= 0$ if and only if
$\lclcoh{i}{B}{M}_{\mbf{d}}\otimes _{R} R_{\mfr{p}}(T)= 0$ for all
$\mfr{p}\in \Spec (R)$. But the formation of local cohomology
commutes with flat base-change:
$\lclcoh{i}{B}{M}_{\mbf{d}}\otimes_{R} R_{\mfr{p}}(T) =
\lclcoh{i}{B_{\mfr{p}}(T)}{M_{\mfr{p}}(T)}_{\mbf{d}}$. Thus
problem 1. above is reduced to a set of similar problems over each
$R_{\mfr{p}}(T)[x_1, ..., x_n ]$, and  $R_{\mfr{p}}(T)$ is a local
ring with infinite residue field.

As for problem 2. the issue is to know whether the $R$-submodule
$P$ of some $Q = M_{\mbf{d}}$ generated by elements from
$S_{\mbf{d}-\mbf{e}}M_{\mbf{e}}$ for all the $\mbf{e}$ in some
given region of $G$, is all of $Q = M_{\mbf{d}}$. We check this by
all localizations $R_{\mfr{p}}(T)$.  We are thus reduced to handling
problem 1 and 2 when $R$ is a local ring with infinite residue
field.

Therefore we must verify a. above.  The one serious place where
the hypothesis ``$K$ is an infinite field'' is used is to assert
the existence of generic forms $g\in S_{\mbf{d}}$ with good
properties, which are then utilized in induction arguments.
Actually Maclagan and Smith only assume that $K$ is a field and
base extend to an infinite field. Ooishi also observed that the
existence of these forms can be shown if the ground ring is a
local ring with infinite residue field.

First we recall some definitions.
\begin{Def}
\label{D:btors} A module $M$ is {\it $B$-torsion} if $M=H^0_B(M)$.
If $M$ is $B$-torsion, then $H^i_B(M)=0$ for all $i>0$.
\end{Def}

\begin{Def}
\label{D:azdv} For an element $g\in S$, we set
\[(0:_M g) = \{ f\in M: gf=0\}.\]
This submodule is zero when $g$ is a non-zero-divisor on $M$.  We
say that $g\in S$ is {\it almost a non-zero-divisor on $M$ } if
$(0:_M g)$ is a $B$-torsion module.
\end{Def}

\begin{Def}
\label{D:umring} We say that a commutative ring $R$ is a {\it
$u$-ring} provided $R$ has the property that an ideal contained in
a finite union of ideals must be contained in one of those ideals;
and a {\it $um$-ring} is a ring $R$ with the property that an
$R$-module which is equal to a finite union of submodules must be
equal to one of them.
\end{Def}

\begin{thm}
\label{T:LSU} $R$  is a $um$-ring if and only if the residue field
$R/{\mfr{p}}$ is infinite for each maximal ideal $\mfr{p}$ of $R$,
and $R$ is a $u$-ring if and only if for each maximal ideal
$\mfr{p}$ of $R$ either the residue field $R/{\mfr{p}}$ is
infinite or the quotient ring $R_{\mfr{p}}$ is a Be\'{z}out ring.
\end{thm}

\begin{proof}
For details see [Theorem 2.2, Theorem 2.3 and Theorem 2.6
\cite{QB}].
\end{proof}

The existence of forms with good properties for induction is the
following (see \cite[Prop. 3.1]{MS}).

\begin{pro}
\label{P:forms} Let $S = R[x_1, ..., x_n]$ be a $G$-graded
polynomial ring where $R$ is a Noetherian local ring with an
infinite residue field. Let $B \subset S$ be a nonzero monomial
ideal attached to a maximal cell $\Gamma$ as explained in
\cite[section 1]{MS} (see also section (\ref{S:intro}) of this
paper). Let $M$ be a finitely generated $G$-graded $S$-module. If
$\mbf{p} \in \mc{K}$, with $\mc{K}$ the Kaehler cone, and $g\in
S_\mbf{p}$ is a sufficiently general form, then $g$ is almost a
non-zero-divisor on $M$.
\end{pro}

\begin{proof}
Since $M'=(0:_M g)$ is $B$-torsion  if each element in $M'$ is
annihilated by some power of $B$. This means $M'_{\mfr{p}}=0$ for
all prime ideals ${\mfr{p}} \subset S$ and $B\nsubseteq
{\mfr{p}}$. That is, $g$ is a non-zero-divisor on the localization
$M_{\mfr{p}}$. Therefore, we only need to show that $g$ is not
contained in any of the associated primes of $M$ except possibly
those which contain $B$. Let
\[
\Ass_+(M)=\{\mfr{p} \in \Ass(M) : \mfr{p} \nsupseteq B\}.
\]
When $n=0$, we have $B = (0)$, $\Ass_+(M)=\emptyset$, and the
assertion is trivial. Suppose $n >0$, and $\Ass_+(M)=\{ \mfr{p}_1,
\dots, \mfr{p}_q\}$ where $\mfr{p}_i$ are $G$-graded ideals.
Suppose  we had an equality of $R$-modules:
\[
S_{\mbf{p}}= \max(R)S_{\mbf{p}} \bigcup (\mfr{p}_1 \cap
S_{\mbf{p}}) \bigcup \dots \bigcup ( \mfr{p}_q \cap S_{\mbf{p}})\]
Because $R$ is an $um$-ring by theorem \ref{T:LSU}, $S_{\mbf{p}}$
would have to be equal to one of the terms in the union.  If
$S_{\mbf{p}}= \max(R)S_{\mbf{p}}$, then by Nakayama's Lemma we
must have $S_{\mbf{p}}=0$. This is not the case: since $\mbf{p}$
is in the Kaehler cone, we have by \cite[Lemma 2.4]{MS} that
\[
B \subseteq \sqrt{\langle S_{\mbf{p}}} \rangle
\]
which is nonzero since $B$ is. Note that Lemma 2.4 of
Maclagan-Smith's paper is valid here - the result is really about
the combinatorics of monomials and has nothing to do with the
ground ring $R$. If say $S_{\mbf{p}}=\mfr{p}_1 \cap S_{\mbf{p}}$,
we would have $\langle S_{\mbf{p}} \rangle \subseteq \mfr{p}_1$,
and hence
\[B \subseteq
\sqrt{\langle S_{\mbf{p}}} \rangle \subseteq \sqrt{\mfr{p}_1}
\subseteq \mfr{p}_1, \] which contradicts the hypothesis on
$\mfr{p}_1$. Any element $g$ of
\[
S_{\mbf{p}} - \max(R)S_{\mbf{p}} \bigcup (\mfr{p}_1 \cap
S_{\mbf{p}}) \bigcup \dots \bigcup ( \mfr{p}_q \cap S_{\mbf{p}})
\]
is almost a non-zero-divisor on $M$.
\end{proof}

\begin{rem}
\label{R:application} Our main application of this section is to
allow coarsening vector  $\mbf{v}$ as in Section \ref{S:vdot},
where $\mbf{v}$ is not a positive coarsening vector. Thus some of
the variables $x_i$ might have $ \mbf{v}\cdot \deg (x_i)= 0$. Our
method is to write $S=R[x_j]$ where $R=K[x_i]$ where $x_i$ are the
variables with $\mbf{v}\cdot \deg(x_i)=0$, and $x_j$ are the
variables with  $ \mbf{v}\cdot \deg (x_j)> 0$, so that in effect
the variables of  $\mbf{v}$-degree zero are treated as constants.
The results of \cite{SVW} refer to vanishing of local cohomology
with respect to the ideal $\mfr{m}$ generated by all variables of
$S$. If we allow some of the variables to have degree zero, then
we must take  $\mfr{m} = \langle x_j\rangle$, the ideal generated
by the variables with strictly positive  $\mbf{v}$-degree. This is
understood when applying the results of \cite{SVW} to this case.
\end{rem}

\section{Examples}\label{S:examples}

\begin{exa}
\label{Ex:classical} Let $S = K[x_1 , ..., x_n ]$ be a polynomial
algebra over a field with standard grading, $\deg (x_i ) = 1$, $B
= \langle x_1 , ..., x_n\rangle$ and $\mc{C} = \{\mbf{1} \}$.
 Let $m \ge 0$ be an integer and set
\[
J = \{J_0 = m, J_1 = m+1, ... \}.
\]
Then an easy calculation using Theorem \ref{P:rsimple} shows that
$\reg _{B, \mc{C}}(J) = m + \N$. On the other hand let $\mc{D} = m
+ \N $, then Theorem \ref{P:drsimple} shows that $\dreg _{B,
\mc{C}}(\mc{D}) = J$ for the set $J$ above. These are exactly the
regularity and degree regions that correspond to each other as in
Theorem \ref{T:BM}.
\end{exa}
The next example is done in \cite[Lemma 3.5]{SVW}. We re-derive it
using Theorem \ref{P:rsimple}.
\begin{exa}
\label{Ex:weighted} Let $S = K[x_1 , ..., x_n ]$ be a polynomial
algebra over a field with a weighted grading, $\deg (x_i ) = a_i >
0$, $B = \langle x_1 , ..., x_n\rangle$ and $c$ is the least
common multiple  of the $a_i$.  The Kaehler cone is $\mc{K} = \N
c$. Take $\mc{C}= \mc{K}$. It is known that (see \cite[Remark
3.3]{MS}) that $\lclcoh{i}{B}{S}_{} = 0$ unless $i = n$. Moreover,
a computation based on \cite[Proposition 3.2]{MS} shows that
$H^n_{B}(S)_w=0$ when $w > -\sum_{i=1}^n a_i$. Therefore
\begin{eqnarray*}
\reg(S) & = & \{u \in \Z: H^i_{B}(S)_w=0,\
\forall i \geq 0,  \  w \in u + \N c[1-i] \}\\
&=& \{ u\in \Z: H^n_{B}(S)_w=0,\  \forall  w \in u+\N c[1-n]
\}\\
&=& \{ u\in \Z: H^n_{B}(S)_w=0,\  \forall  w \in u+(1-n)c+\N c \},
\end{eqnarray*}
which shows that
\[
\{ u \in \Z: u \geq (n-1)c - \sum_{i=1}^n a_i +1 \} \subseteq
\reg(S).
\]
If $M$ is a graded $S$-module with a minimal resolution of type
$J$, then $\reg_{B, \mc{C}} (M)$ contains $\reg_{B, \mc{C}} (J)$
by Lemma \ref{C:regJ} and this last one is seen by Theorem
\ref{P:rsimple} to be the set of integers that simultaneously
belong to all the sets
\[
d_p + \mr{reg}(S) -pc + \N c, \ \ \mr{for\ all\ } d_p \in J_p.
\]
Using the information on $\mr{reg}(S)$ above, we see that
$\mr{reg}(M)$ contains the interval
\[
\{m \in \Z : m \ge d_p
+ (n - p -1)c - \sum a_i + 1 \text{\ for \ all\ }
0 \le p \le s, \ d_p \in J_p \}.
\]
\end{exa}

\begin{exa}\label{E:multiproj}
Let $S=\C[x_{11}, \dots, x_{1m_1}; \dots; x_{l1}, \dots
x_{lm_l}]$, with $\deg(x_{ij})=\mbf{e}_i \in \Z^l=G$, where $\{
\mbf{e}_i\}$ is the standard basis of $\R ^ l$. This is the
homogeneous coordinate ring of $\proj{m_1 -1}\times ...\times
\proj{m_l -1}$. The ideal $B$ is generated by all $x_{1,
i_{1}}...x_{l, i_{l}}$. Set $\mc{C} = \{ \mbf{e}_i\}$. The theory
developed in Section \ref{S:family} applies, with $B_i = \langle
x_{i1}, \dots, x_{im_i}\rangle $. For each  $\emptyset \neq I
\subseteq \{1, 2, \dots, l\}$, we define $B_I=\sum _{i \in I
}B_i$. Define a vector grading for this family via
$\mbf{v}_I=\sum_{i\in I} \mbf{e}_i$. We have
\[H^i_{B_I}(S)_{\mbf{d}}=0, \text{ for \ all  \  }\mbf{d}\cdot \mbf{v}_I
\geq 1-i, \text{ for each } I.\] This implies
\[ \reg_{B_*, \mc{C}, \mbf{v}_*}(S)=\N^l.\]
Let $M$ be a $\Z^l$-graded $S$-module with minimal resolution
\[
\small{
\begin{CD}
0 @>>> \displaystyle{\bigoplus_{\mbf{d}_s \in J_s}}
S(-\mbf{d}_{s}) @>>> \cdots @>>>
\displaystyle{\bigoplus_{\mbf{d}_0 \in J_0}} S(-\mbf{d}_0) @>>> M
@>>> 0
\end{CD}}
\]
If $\mbf{p} \in \reg_{B_*,\mc{C},  \mbf{v}_*}(M)$ then the degrees
of the syzygies satisfy the bounds $\mbf{d}_{i}\cdot \mbf{v}_I
\leq \mbf{p} \cdot \mbf{v}_I +i$ for all $\emptyset \neq I =\{ 1,
2, \dots, l\}$. (compare \cite[Section 4]{HW}). This can be seen
as an application of Proposition \ref{P:vres}. First note that in
the Definition $\ref{D:Bv}$ we have $m(1, i) = 1$ for all $i$.
Therefore if $\mbf{p} \in \reg_{B_*,\mc{C},  \mbf{v}_*}(M)$ the
condition
\[
\lclcoh{i}{B_{I}}{M}_{\mbf{d}}=0\ \text{for\ all\ } \deg
_{\mbf{v}_{I}} (\mbf{d})\ge\deg _{\mbf{v}_{I}} (\mbf{m})
                 +(1 - i)
\]
shows that $b_I := \mbf{v}_I \cdot \mbf{p}\ge \vIregnum (M)$ since
$\vIregnum (M)$ is the least integer which gives the above
vanishing statement for all $i \ge 0$. Since $c _{\mbf{v}_{I}} =
1$, $\deg _{\mbf{v}_I} (\mbf{p}) = \mbf{v}_I \cdot \mbf{p}$
satisfies the hypotheses of Proposition \ref{P:vres}. From that
proposition and Corollary \ref{C:resnum} we conclude that the
$i$th syzygies are in $K_i (\mbf{v}_* ,  b_*)$.

Note that $s_{\mbf{v}_{I}} = 1$, so that this region is just
$\mbf{d}_{i}\cdot \mbf{v}_I \leq \mbf{p} \cdot \mbf{v}_I +i$.
\end{exa}

\begin{exa}
\label{Ex:hirz} Let $S=\C[x_1, x_2, x_3, x_4]$ be the coordinate
ring of the Hirzebruch surface $\F_t$ with $t\geq 0$, $G=\Z^2$,
$\deg(x_1)=(1, 0)$, $\deg(x_2)=(-t, 1)$, $\deg(x_3)=(1,0)$, and
$\deg(x_4)=(0,1)$. $B=(x_1x_2, x_1x_4, x_2 x_3, x_3x_4)$. Let
$\mc{C} = \{(1, 0), (0, 1) \}$. Then (\cite[Example 4.3]{MS})
\[
\reg (S) = \begin{cases}
           \N ^2 \ \ \ \ \ \ \ \ \ \ \ \ \
\ \ \ \ \ \ \ \ \ \ \ \ \ \ \ \ \ \ \ \ \ \ \ \ \ \
\mr{for\ } t = 0, 1;\\

 ((t-1, 0)+ \N ^2)\cup ((0, 1)+\N ^2) \ \   \mr{for \ } t \ge 2.
           \end{cases}
\]
Natural choices of vector coarsenings are $\mbf{v} = (0, 1)$ and
$\mbf{w} = (1, t)$. The $\mbf{v}$-gradings behave just like
Example \ref{Ex:classical} but with the variables $x_1$ and $x_3$
assigned degree $0$.  The $\mbf{w}$-gradings behave just like
Example \ref{Ex:weighted} the variable $x_2$ assigned degree $0$,
the variables $x_1$ and $x_3$ given degree $1$ and the variable
$x_4$ given degree $t$. By imposing regularity conditions with
respect to these vector gradings we can describe syzygy regions by
inequalities for some $a,\  b, \  c, \  d$:
\begin{eqnarray*}
\mbf{d}_i\cdot \mbf{v}  & \leq & \mbf{p} \cdot \mbf{v} +a i + b, \\
\mbf{d}_{i}\cdot \mbf{w} & \leq & \mbf{p} \cdot \mbf{w} +c i + d.
\end{eqnarray*}
Let
 $Z=\V(B)=Z_1 \cup Z_2=\V(B_1) \cup \V(B_2)$, where
$B_1=(x_1, x_3)$, $B_2=(x_2, x_4)$. This is the Batryev
decomposition in this case. However, there is not an orthogonal
choice of vector gradings for the ideals $B_1$, $B_2$ in the sense
of Section \ref{S:family} so it does not seem that we can apply
the results of that section here.
\end{exa}

\end{document}